%% file: obstructions_for_extension.tex
\documentclass[envcountsect,envcountsame]{llncs}
\usepackage{makeidx}
\usepackage{amsfonts}
\usepackage{amsmath}
\usepackage{amssymb}
\usepackage{amscd}
\usepackage[dvips]{graphicx}
\usepackage{algorithmic}
\usepackage{algorithm}
\usepackage{comment}
\usepackage{color}
\usepackage{array,amssymb,amsmath,amsbsy}
\usepackage{enumerate}
\usepackage[bottom]{footmisc}
\usepackage[all,cmtip]{xy}
\usepackage{amsbsy}
\usepackage{a4wide}
\usepackage{hyperref}
\usepackage{simplemargins}

\setallmargins{1in}

\makeatletter
\renewcommand*\l@author[2]{}
\renewcommand*\l@title[2]{}
\makeatletter
\newcommand{\nocontentsline}[3]{}
\newcommand{\tocless}[2]{\bgroup\let\addcontentsline=\nocontentsline#1{#2}\egroup}

\newenvironment{packed_enum}{
	\begin{enumerate}
		\setlength{\itemsep}{1pt}
	    \setlength{\parskip}{0pt}
		\setlength{\parsep}{0pt}
}{\end{enumerate}}

\newenvironment{packed_itemize}{
	\begin{itemize}
		\setlength{\itemsep}{1pt}
	    \setlength{\parskip}{0pt}
		\setlength{\parsep}{0pt}
}{\end{itemize}}

\newcommand{\heading}[1]{\medskip\par\noindent{\bf #1}}

\def\computationproblem#1#2#3#4{
	\vskip 1ex
	\begin{center}
	\fbox{\begin{tabular}{rp{#4}}
	{\bf Problem:\enspace}&#1\\
	{\bf Input:\enspace}&#2\\
	{\bf Output:\enspace}&#3\\
	\end{tabular}}
	\end{center}
	\vskip 1ex
}

 \def\calR{{\cal R}} \def\calS{{\cal S}}

\def\cNP{\hbox{\rm \sffamily NP}}

\def\O{\mathcal{O}{}}

\def\wlt{\vartriangleleft}
\def\wgt{\vartriangleright}

\def\clql{\mathord\curvearrowleft}
\def\clqr{\mathord\curvearrowright}
\def\da{\downarrow}
\def\cp{\textrm{\rm cp}}

\def\inter#1{\left<#1\right>}

\def\endlft{{\mathord{\reflectbox{\ensuremath{\scriptstyle\mapsto}}}}}
\def\startrt{{\mathord{\mapsto}}}
\def\lft{\leftarrow}
\def\rt{\rightarrow}

\def\ext{\textsc{RepExt}}

\def\int{\hbox{\bf \rm \sffamily INT}}
\def\pint{\hbox{\bf \rm \sffamily PROPER INT}}
\def\uint{\hbox{\bf \rm \sffamily UNIT INT}}
\def\circle{\hbox{\bf \rm \sffamily CIRCLE}}
\def\fun{\hbox{\bf \rm \sffamily FUN}}
\def\perm{\hbox{\bf \rm \sffamily PERM}}

\def\chor{\hbox{\bf \rm \sffamily CHOR}}
\def\atfree{\hbox{\bf \rm \sffamily AT-FREE}}

\newcounter{lth}
\title{Minimal Obstructions for Partial Representations\\of Interval Graphs\thanks{The conference
version of this paper appeared in ISAAC 2014~\cite{ks}. For an interactive structural diagram of
this paper, see \url{http://pavel.klavik.cz/orgpad/minobstr.html} (supported for Firefox and Google
Chrome).}}

\author{Pavel Klav\'{\i}k\inst{1}${}^,$\thanks{%
			Supported by CE-ITI (P202/12/G061 of GA\v{C}R) and Charles University as GAUK 196213.
			A part of this work was carried out during a visit at ULB in Brussels, funded by a
			EUROCORES Short Term Visit grant.}
		\and Maria Saumell\inst{2}${}^,$\thanks{%
			Supported by project LO1506 of the Czech Ministry of Education, Youth and Sports, project NEXLIZ - CZ.1.07/2.3.00/30.0038 co-financed by
			European Social Fund and the state budget of Czech Republic,
			and ESF EuroGIGA project ComPoSe as F.R.S.-FNRS - EUROGIGA NR 13604.}
		}

\institute{Computer Science Institute, Charles University in Prague,\\Czech Republic.
			E-mail: \texttt{klavik@iuuk.mff.cuni.cz}\\[0.75ex]\and
		Department of Mathematics and European Centre of Excellence NTIS,\\University of West Bohemia, Pilsen, Czech Republic.
			E-mail: \texttt{saumell@kma.zcu.cz}.}

\begin{document}
\maketitle

\begin{abstract}
\emph{Interval graphs} are intersection graphs of closed intervals. A~generalization of recognition
called \emph{partial representation extension} was introduced recently. The input gives an interval
graph with a \emph{partial representation} specifying some pre-drawn intervals.  We ask whether the
remaining intervals can be added to create an \emph{extending representation}.  Two linear-time
algorithms are known for solving this problem.

\hskip 1em In this paper, we characterize the \emph{minimal obstructions} which make partial
representations non-extendible. This generalizes Lekkerkerker and Boland's characterization of
the minimal forbidden induced subgraphs of interval graphs. Each minimal obstruction consists of a
forbidden induced subgraph together with at most four pre-drawn intervals. A Helly-type result
follows: A partial representation is extendible if and only if every quadruple of pre-drawn
intervals is extendible by itself. Our characterization leads to a linear-time certifying algorithm for
partial representation extension.
\end{abstract}

\tableofcontents
\newpage

\input introduction
\input definition_of_obstructions
\input maximal_cliques

\input interval_orders
\input strategy
\input leaves

\input p_nodes

\input q_nodes
\input main_result
\input conclusions

\heading{Acknowledgment.} We would like to thank to Jaroslav Ne\v{s}et\v{r}il for suggesting us the
problem of minimal obstructions for partial representation extension.

\bibliographystyle{plain}
\bibliography{obstructions_for_extension}

\end{document}

%% file: introduction.tex
\section{Introduction} \label{sec:intro}

The main motivation for graph drawing and geometric representations is finding ways
to visualize some given data efficiently. The most famous representations are plane drawings, in
which we draw a graph in the plane and we want to avoid (or minimize) crossings of edges.
However, for certain types of graphs, \emph{intersection representations} are more suitable.  They
represent each vertex by a geometrical object and encode the edges by intersections.

\subsection{Interval Graphs}

\begin{figure}[b!]
\centering
\includegraphics{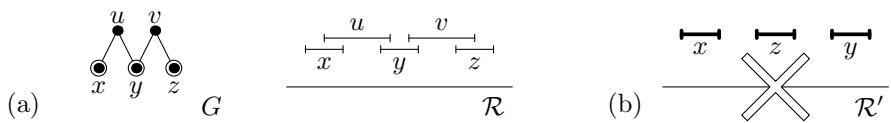}
\caption{(a) An interval graph $G$ with one of its interval representations $\calR$.  (b) A partial
representation $\calR'$ with pre-drawn intervals $\inter x'$, $\inter y'$ and $\inter z'$. It is
non-extendible since $\inter u$ cannot be placed. In all figures, we depict pre-drawn intervals in
bold.}
\label{fig:int_example}
\end{figure}

The most studied class of intersection graphs are \emph{interval graphs} (\int), defined by
H\'ajos~\cite{hajos_interval_graphs} in 1957. An \emph{interval representation} $\calR$ is a
collection of closed intervals $\bigl\{\inter x : x \in V(G)\bigr\}$ where $\inter x \cap \inter y
\ne \emptyset$ if and only if $xy \in E(G)$. A graph is an interval graph if it has an interval
representation; see Fig.~\ref{fig:int_example}a.

Interval graphs have many applications. Already in 1959, Benzer~\cite{benzer_interval_graphs} used
them in his experimental study of DNA. For some time, interval graphs played an important role for
the DNA hybridization~\cite{karp_hybridization}, in which short pieces of DNA are studied
independently. Further applications include scheduling, psychology, archaeology,
etc.~\cite{stoffers,roberts_discrete_models,kendall}.

Interval graphs also have nice theoretical properties. They are perfect and closely related to
path-width decompositions. They can be recognized in linear
time~\cite{PQ_trees,LBFS_int,finding_lb_graphs}, and many hard combinatorial problems are
polynomially solvable for interval graphs.  Fulkerson and Gross~\cite{maximal_cliques} characterized
them by consecutive orderings of maximal cliques (see Section~\ref{sec:maximal_cliques} for
details).  This lead Booth and Lueker~\cite{PQ_trees} to the construction of PQ-trees, which are an
efficient data structure to deal with consecutive orderings, and have many other applications. 

\begin{figure}[t!]
\centering
\includegraphics{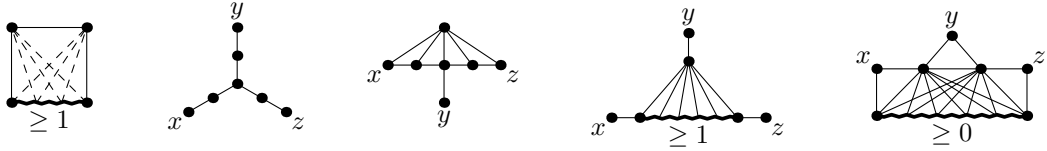}
\caption{Five types of LB obstructions which are the minimal forbidden induced subgraphs of \int.
The bold curly lines correspond to induced paths with denoted minimal lengths. The leftmost
obstructions are induced cycles of length four or more.  The remaining four types of obstructions
are minimal asteroidal triples $(x,y,z)$ which are chordal graphs.}
\label{fig:LB_obstructions}
\end{figure}

\emph{Chordal graphs} (\chor) are graphs with no induced cycle of length four or more, alternatively
intersection graphs of subtrees of trees. Three vertices form an \emph{asteroidal triple} if there
exists a path between every pair of them avoiding the neighborhood of the third vertex.
\emph{Asteroidal triple-free graphs} (\atfree) are graphs containing no asteroidal triples.
Lekkerkerker and Boland~\cite{lb_graphs} characterized interval graphs as $\int = \chor \cap
\atfree$. They described this characterization by the minimal forbidden induced subgraphs given in
Fig.~\ref{fig:LB_obstructions} which we call \emph{Lekkerkerker-Boland obstructions} (LB).

\subsection{Partial Representation Extension}

The partial representation extension problem was introduced by Klav\'{\i}k et al.~\cite{kkv}.  A
\emph{partial representation} $\calR'$ of $G$ is an interval representation $\bigl\{\inter x' : x \in
V(G')\bigr\}$ of an induced subgraph $G'$ of $G$. The vertices of $G'$ and the intervals of $\calR'$
are called \emph{pre-drawn}.  A representation $\calR$ of $G$ \emph{extends} $\calR'$ if and only if
it assigns the same intervals to the vertices of $G'$, i.e., $\inter x = \inter x'$ for every $x \in
V(G')$.  For an example, see Fig.~\ref{fig:int_example}b.

\computationproblem
{Partial Representation Extension -- $\ext(\int)$}
{A graph $G$ and a partial representation $\calR'$ of $G'$.}
{Is there an interval representation of $G$ extending $\calR'$?}
{9.25cm}

\noindent The first polynomial-time algorithm, running in $\O(n^2)$ time, was given in~\cite{kkv}.
Currently, there are two different linear-time algorithms~\cite{blas_rutter,kkosv} for this problem.
	
We note that the partial representation extension problems have been considered also for other classes
of intersection graphs. A linear-time algorithm for proper interval graphs and an almost quadratic-time
algorithm for unit interval graphs are given in~\cite{kkorssv}, and improved to quadratic time
in~\cite{soulignac}. The partial representation extension problems are polynomial-time solvable for
$k$-nested interval graphs (classes generalizing proper interval graphs), but \cNP-hard for
$k$-length interval graphs (classes generalizing unit interval graphs), even for $k=2$~\cite{kos}.
Polynomial-time algorithms are further known for circle graphs~\cite{cfk}, and
permutation and function graphs~\cite{kkkw}. The partial representation extension problems for
chordal graphs~\cite{kkos} and contact representations of planar graphs~\cite{contact_planar_ext}
are \cNP-hard. Notable graph classes for which the complexity of the partial representation
extension problem is open are circular-arc graphs and trapezoid graphs.

Outside intersection graphs, the similar problem was considered even sooner for planar
graphs.  Partially embedded planar graphs can be extended in linear time~\cite{extending_planar}.
Even though every planar graph has a straight-line embedding, extension of such embeddings is
\cNP-hard~\cite{straight_planar_npcomplete}. Kuratowski's characterization of minimal forbidden
minors was extended to partially embedded planar graphs by Jel\'{\i}nek et
al.~\cite{kuratowski_ext}. Our research has a similar spirit as this last result.

\begin{figure}[b!]
\centering
\includegraphics{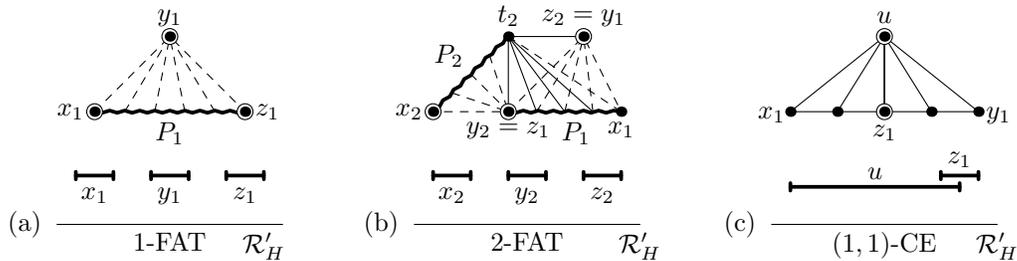}
\caption{Three examples of minimal obstructions, each consisting of a graph $H$ and a non-extendible
partial representation $\calR'_H$. Curly lines denote induced paths and dashed edges are non-edges.
The obstructions (a) and (b) are the first two $k$-FAT obstructions, and (c) is the simplest
$(k,\ell)$-CE obstruction.}
\label{fig:minimal_obstructions}
\end{figure}

\subsection{Our Results}

In this paper, we generalize the characterization of Lekkerkerker and Boland~\cite{lb_graphs} to
describe minimal obstruction which make partial representations non-extendible.  Each obstruction
consists of a small graph and its non-extendible partial representation.  Aside LB obstructions, we have
two trivial obstructions, called SE, and ten infinite classes of minimal obstructions. The main
class, called $k$-FAT obstructions, has three wrongly ordered disjoint pre-drawn intervals
$\inter{x_k}'$, $\inter{y_k}'$, and $\inter{z_k}'$. The obstruction consists of a zig-zag structure
with $k$ levels where the last level cannot be placed.  See Fig.~\ref{fig:minimal_obstructions}a and
b for $1$-FAT and $2$-FAT obstructions. There are eight other infinite classes derived from $k$-FAT
obstructions by adding a few vertices and having different vertices pre-drawn. The last infinite
class of $(k,\ell)$-CE obstructions consists of a $k$-FAT obstruction glued with an $\ell$-FAT
obstruction and contains only two pre-drawn vertices; see Fig.~\ref{fig:minimal_obstructions}c for a
$(1,1)$-CE obstruction.  We formally define these minimal obstructions in
Section~\ref{sec:definition_of_obstructions}.

\begin{theorem} \label{thm:repext_obstructions}
A partial representation $\calR'$ of $G$ is extendible if and only if $G$ and $\calR'$ contain no
LB, SE, $k$-FAT, $k$-BI, $k$-FS, $k$-EFS, $k$-FB, $k$-EFB, $k$-FDS, $k$-EFDS, $k$-FNS and
$(k,\ell)$-CE obstructions.
\end{theorem}

Since every minimal obstruction contains at most four pre-drawn intervals, we get the following
Helly-type result as a straightforward corollary:

\begin{corollary} \label{cor:helly_obstructions}
A partial representation is extendible if and only if every quadruple of pre-drawn intervals is
extendible by itself.
\end{corollary}

All known algorithms for the partial representation extension problems
\cite{kkv,kkkw,blas_rutter,kkosv,kkorssv,kkos} are able to certify solvable instances by outputting
an extending representation. Using our minimal obstructions, we construct the first algorithm for
partial representation extension certifying also non-extendible partial
representations.\footnote{Formally speaking, a polynomial-time algorithm certifies
unsolvable instances by outputting ``no'' and by a proof of its correctness. Our algorithm
outputs a simple proof that a given partial representation is non-extendible in terms of a
minimal obstruction. This proof can be independently verified which is desirable.}

\begin{corollary} \label{cor:certifying_algorithm}
Assume that the input gives the endpoints in a partial representation sorted from left to right.
Then there exists an $\O(n+m)$ certifying algorithm for the partial representation extension problem,
where $n$ is the number of vertices and $m$ is the number of edges of the input graph. If the answer
is ``yes'', it outputs an extending representation. If the answer is ``no'', it detects one of the
minimal obstructions.
\end{corollary}

\heading{Outline.} In Section~\ref{sec:definition_of_obstructions}, we define the minimal
obstructions which make a partial representation non-extendible. In
Section~\ref{sec:maximal_cliques}, we introduce the standard tools for working with interval graphs:
the characterization of Fulkerson and Gross by linear orderings of maximal cliques, and the related
data structure called an MPQ-trees which stores all feasible orderings.

In Section~\ref{sec:interval_orders}, we restate the characterization of~\cite{kkosv}: a partial
representation is extendible if and only if there exists a feasible ordering of the maximal cliques
which extends a certain partial ordering $\wlt$. Therefore, to solve $\ext(\int)$, we test whether
the MPQ-tree can be reordered according $\wlt$.

Based on this, in Section~\ref{sec:strategy}, we build our strategy for showing that every
non-extendible partial representation contains one of the defined minimal obstructions. If
reordering of the MPQ-tree according to $\wlt$ fails, then it fails in a leaf, in a P-node, or a
Q-node. We deal with these three cases in Sections~\ref{sec:leaves}, \ref{sec:p_nodes},
and~\ref{sec:q_nodes}, respectively, where the last one is most involved. In
Section~\ref{sec:main_result}, we put these results together and establish
Theorem~\ref{thm:repext_obstructions} and Corollaries~\ref{cor:helly_obstructions}
and~\ref{cor:certifying_algorithm}.

We conclude with a discussion of our results and several open problems.

\subsection{Preliminaries}

For a graph $G$, we denote by $V(G)$ its vertices and by $E(G)$ its edges. We denote the
\emph{closed neighborhood} of $x$ by $N[x]$. Maximal cliques are denoted by the letters $a$ to $f$,
and vertices by the remaining letters. For $A \subseteq V(G)$, we denote by $G[A]$ the subgraph
induced by $A$. Similarly, for $A \subseteq V(G')$, we denote by $\calR'[A]$ the partial
representations which only contains the pre-drawn intervals in $A$. By $P_{x,y}$ we denote an
induced path from $x$ to $y$; its length is the number of edges.

For an interval $\inter x$, we denote its left endpoint by $\ell(x)$ and its right endpoint by
$r(x)$. If $r(x) < \ell(y)$, we say that $\inter x$ is \emph{on the left} of $\inter y$ and $\inter
y$ is \emph{on the right} of $\inter x$. We say that $\inter y$ is \emph{between} $\inter x$ and
$\inter z$ if $\inter x$ is on the left of $\inter y$ and $\inter z$ is on the right of $\inter y$,
or vice versa. We also work with open intervals, for which the inequalities are non-strict.

We conclude with a list of the remaining notation.  In Section~\ref{sec:maximal_cliques}, we define
MPQ-trees, $s(N)$, $s_i(Q)$, $s^{\lft}_u(Q)$, $s^{\rt}_u(Q)$, $G[T]$, $G[N]$, $T[N]$, and
$Q$-monotone paths. In Section~\ref{sec:interval_orders}, we define $\clql(a)$, $\clqr(a)$, $I_a$,
$\wlt$, the flip operation, $P^{\startrt}(a)$, and $P^{\endlft}(a)$.

%% file: definition_of_obstructions.tex
\section{Definition of Minimal Obstructions} \label{sec:definition_of_obstructions}

In this section, we formally define all twelve classes of minimal obstructions which make a partial
representation non-extendible. 

\heading{Definition.}
Every \emph{obstruction} consists of a graph $H$ and a non-extendible partial representation $\calR'_H$.
This obstruction is \emph{contained} in $G$ and $\calR'$ if (i) $H$
is an induced subgraph of $G$, (ii)~the pre-drawn vertices of $H$ are mapped to pre-drawn vertices of
$G$, and (iii) the endpoints in $\calR'_H$ have the same left-to-right order as the endpoints of the
corresponding pre-drawn vertices in $\calR'$. For instance, the partial representation in
Fig.~\ref{fig:int_example}b contains a $1$-FAT obstruction, given in Fig.~\ref{fig:minimal_obstructions}a.

We give $H$ by \emph{descriptions} using finitely many vertices, edges, and induced paths. For inner
vertices of the induced paths, we specify their adjacencies with the remainder of $H$. Since these
induced paths do not have fixed lengths, each description having at least one induced path defines
an infinite class of forbidden subgraphs $H$. Unlike LB obstructions, most classes of minimal
obstructions need infinitely many different descriptions. For instance, each FAT obstruction has $k$
induced paths, and different values of $k$ need different descriptions.

If $H$ contains an induced path $P_{x,y}$, and $x$ and $y$ are allowed to be adjacent, then
$P_{x,y}$ can be a single edge.  When $N[x] = N[y]$, we allow the length of $P_{x,y}$ to be zero,
i.e., $x=y$.

\heading{Minimality.} An obstruction is \emph{minimal} if $\calR'_H$ becomes extendible when any
vertex is removed or any pre-drawn interval is made \emph{free} by removing it from the partial
representation $\calR'_H$. 

\subsection{List of Minimal Obstructions}

In Fig.~\ref{fig:LB_obstructions}, we have already described minimal LB obstructions
of~\cite{lb_graphs} with $\calR'_H = \emptyset$.  There are eleven other other classes of minimal
obstructions we describe now. 

\begin{figure}[t!]
\centering
\includegraphics{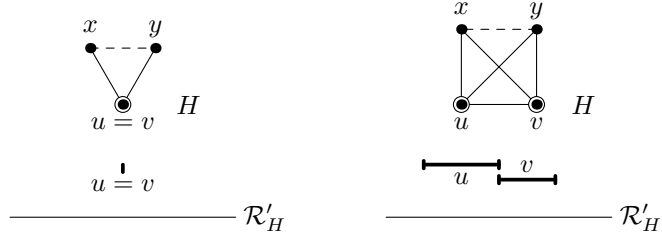}
\caption{The SE obstructions, on the left with $u=v$, on the right with $u \ne v$.}
\label{fig:SE}
\end{figure}

\begin{figure}[b!]
\centering
\includegraphics{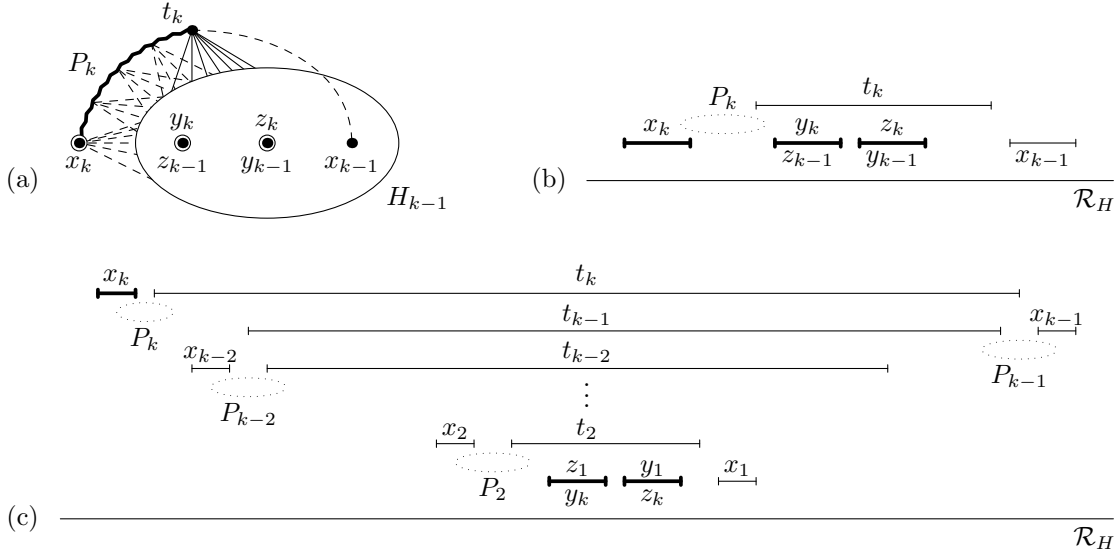}
\caption[ ]{(a) A $k$-FAT obstruction is created from a $(k-1)$-FAT obstruction. It consists of the
vertices $x_1,\dots,x_k$, $t_2,\dots,t_k$, $y_k$, $z_k$, and the induced paths $P_1,\dots,P_k$. The
adjacencies are defined inductively.\\
(b) In every representation $\calR_H$, the pre-drawn interval $\inter{x_k}'$ together with $P_k$ and
$t_k$ forces $\inter{x_{k-1}}$ to be placed on the right of $\inter{z_k}'$. Therefore, the induced
$(k-1)$-FAT obstruction is forced.\\
(c) The global zig-zag pattern forced by a $k$-FAT obstruction, with $k$ nested levels going across
$\inter{y_k}'$ and $\inter{z_k}'$. It is an obstruction since $P_1$ going from $x_1$ to $z_1$ with
all inner vertices non-adjacent to $y_1$ cannot be placed.}
\label{fig:k_FAT}
\end{figure}

\heading{SE obstructions.} We start with two simple \emph{shared endpoint obstructions} which deal
with shared endpoints in $\calR'$; see Fig.~\ref{fig:SE}. We have two pre-drawn vertices $u$ and
$v$ such that $r(u) = \ell(v)$ (possibly $u=v$, so only one interval may be pre-drawn). Further,
there are two non-adjacent vertices $x$ and $y$, both adjacent to $u$ and $v$.  If $u \ne v$, the
minimality requires that $\ell(u) < \ell(v) = r(u) < r(v)$.

\heading{$\boldsymbol{k}$-FAT obstructions.} The class of \emph{forced asteroidal triple obstructions}
is defined inductively; the first two obstructions $1$-FAT and $2$-FAT are depicted in
Fig.~\ref{fig:minimal_obstructions}a and b, respectively.

A $1$-FAT obstruction consists of three pre-drawn non-adjacent vertices $x_1$, $y_1$ and $z_1$ such
that $\inter{y_1}'$ is between $\inter{x_1}'$ and $\inter{z_1}'$. Further, $x_1$ and $z_1$ are
connected by an induced path $P_1$ and $y_1$ is non-adjacent to the inner vertices of $P_1$. See
Fig.~\ref{fig:minimal_obstructions}a.

A $k$-FAT obstruction is defined as follows; see Fig.~\ref{fig:k_FAT}a.  Let $H_{k-1}$ be the graph
of a $(k-1)$-FAT obstruction.  To get $H_k$, we add to $H_{k-1}$ two vertices $x_k$ and $t_k$
connected by an induced path $P_k$.   Concerning edges, $t_k$ is adjacent to all vertices of
$H_{k-1}$, except for $x_{k-1}$. All vertices of $H_{k-1}$ are non-adjacent to $x_k$ and to the
inner vertices of $P_k$.  Further, for $k > 1$, we allow $P_1$ to be a single edge, so $x_1$
can be adjacent to $z_1$.

We put $y_k = z_{k-1}$ and $z_k = y_{k-1}$.  A $k$-FAT obstruction has three pre-drawn vertices
$x_k$, $y_k$ and $z_k$ such that $\inter{y_k}'$ is between $\inter{x_k}'$ and $\inter{z_k}'$.  The
role of $\inter{x_k}'$, $P_k$ and $t_k$ is to force $\inter{x_{k-1}}$ to be placed on the other side of
$\inter{z_k}' = \inter{y_{k-1}}'$ than $\inter{y_k}' = \inter{z_{k-1}}'$, thus forcing the
$(k-1)$-FAT obstruction of $H_{k-1}$; see Fig.~\ref{fig:k_FAT}b. The global structure forced by a
$k$-FAT obstruction is depicted in Fig.~\ref{fig:k_FAT}c. 

\begin{figure}[b!]
\centering
\includegraphics{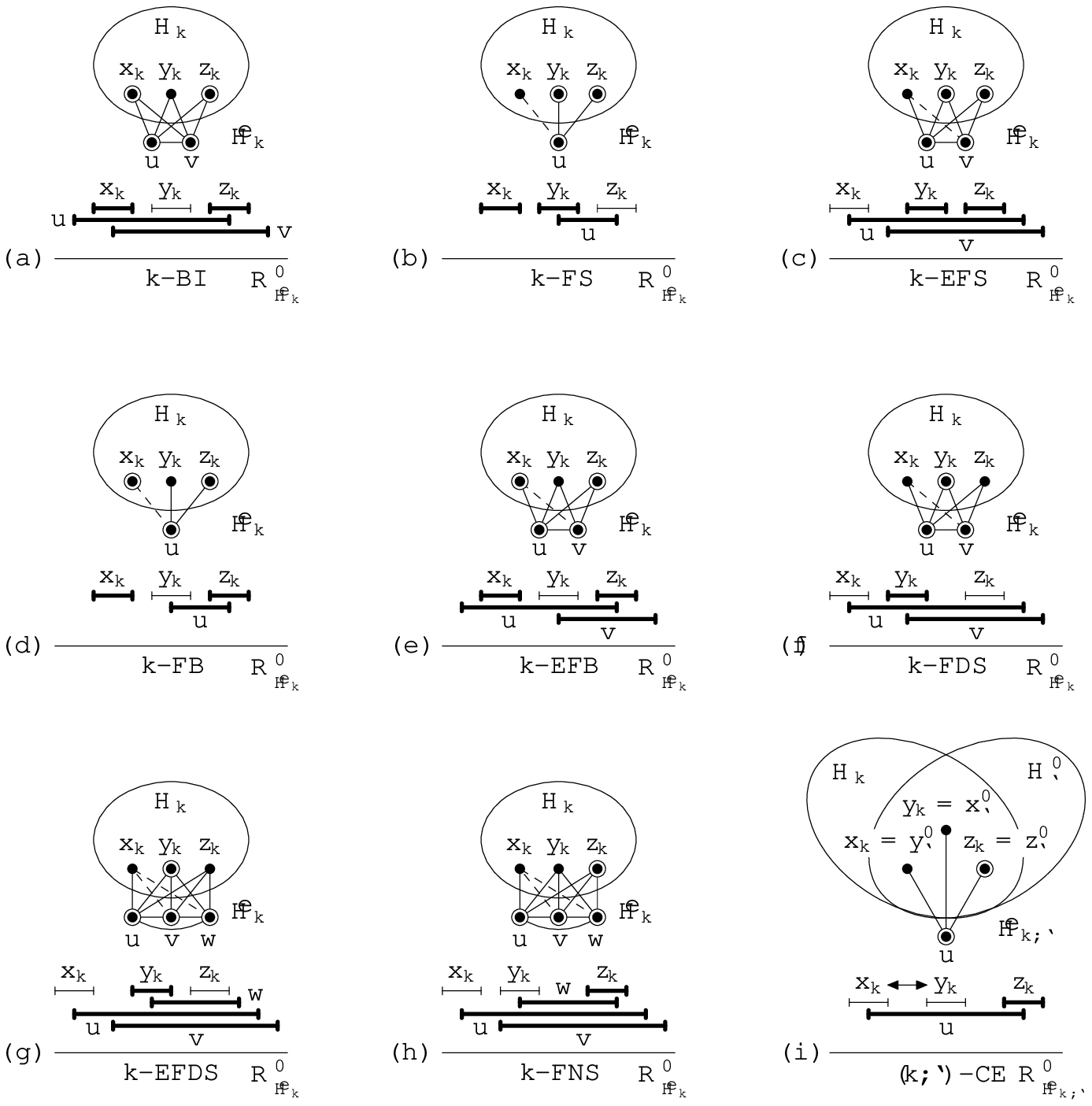}
\caption[ ]{Nine classes of obstructions derived from $k$-FAT obstructions.}
\label{fig:overview_of_obstructions}
\end{figure}

The remaining nine classes, depicted in Fig.~\ref{fig:overview_of_obstructions}, are derived from
$k$-FAT obstructions. Let $H_k$ denote the graph of a $k$-FAT obstruction. With exception of the
last class of $(k,\ell)$-CE obstructions, we create the graphs $\widetilde H_k$ of these
obstructions by adding a few vertices to $H_k$. For $k=1$, when one of $x_1$ and $z_1$ is not
pre-drawn, we also allow $x_1$ to be adjacent to $z_1$. We also consider all of these
obstructions horizontally flipped.

It is possible that the added vertices already belong to $H_k$; for instance, a $k$-BI obstruction
may have $u = t_k$ or $v = t_k$.  Also, we do not specify in details the edges between the added
vertices and $H_k \setminus \{x_k, y_k, z_k\}$.  An accurate description would be too lengthy and
the reader may derive it from Fig.~\ref{fig:k_FAT}c.  We describe all minimal $(k,\ell)$-CE
obstructions in details.

\heading{$\boldsymbol{k}$-BI obstructions.} The class of \emph{blocked intersection
obstructions} is shown in Fig.~\ref{fig:overview_of_obstructions}a. To create $\widetilde H_k$ from $H_k$, we add
two vertices $u$ and $v$ adjacent to $x_k$, $y_k$ and $z_k$. Then the partial
representation contains four pre-drawn vertices $x_k$, $z_k$, $u$ and $v$.
We have $\ell(u) \le \ell(v) < r(u) \le r(v)$, $\inter{x_k}'$ covering $\ell(v)$, and
$\inter{z_k}'$ covering $r(u)$. We allow $u = v$. 

\begin{figure}[t!]
\centering
\includegraphics{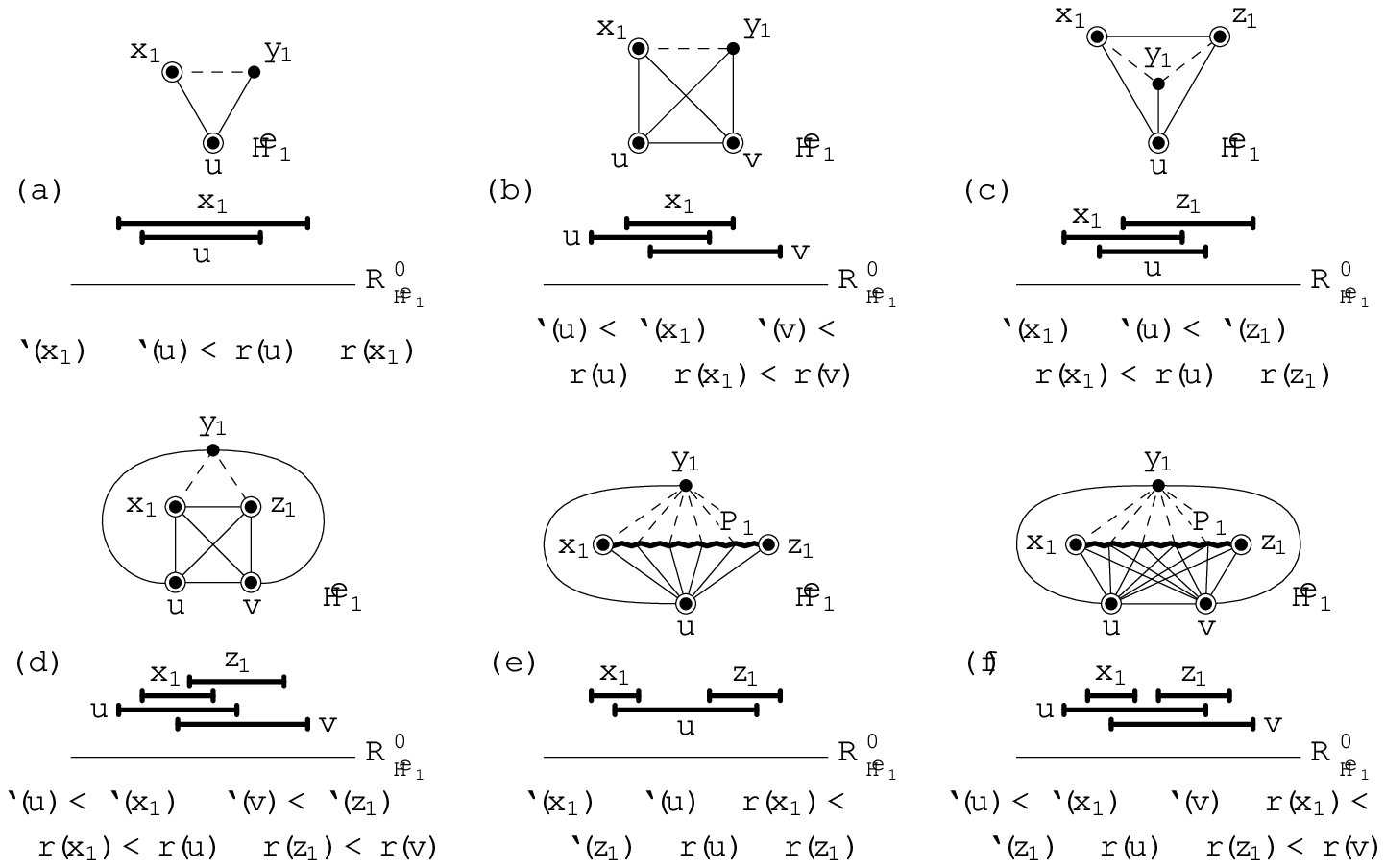}
\caption{All minimal 1-BI obstructions are depicted.  Each gives several $1$-BI obstructions with
different $\calR'_{\widetilde H_1}$, since there are several possible orderings of the endpoints
satisfying the given inequalities.  We have $x_1=z_1$ for (a) and (b), and $x_1 \ne z_1$ otherwise.
We have $u=v$ for (a), (c) and (e), and $u \ne v$ otherwise.}
\label{fig:BI_obstructions}
\end{figure}

The minimality further implies that $k \le 2$. Indeed, a $k$-BI obstruction with $k>3$ contains a smaller
$(k,1)$-CE obstruction by removing $v$ and freeing $\inter{x_k}'$ (this follows from
Lemma~\ref{lem:k_fat_contains_1_fat}). Concerning $1$-BI, we allow $x_1 = z_1$. We illustrate all possible cases only for $1$-BI
obstructions, so the reader can understand the complexity of these classes; see
Fig.~\ref{fig:BI_obstructions}.  For $k=2$, we know by Lemma~\ref{lem:k_fat_contains_1_fat} that
$x_2$ is adjacent to $t_2$. The pre-drawn intervals are as follows:
\begin{align*}
\ell(x_2) \le \ell(u) \le r(x_2)&<\ell(z_2) \le r(u) \le r(z_2), & \text{for $u=v$,}\\
\ell(u) < \ell(x_2) \le \ell(v) \le r(x_2)&<\ell(z_2) \le r(u) \le r(z_2) < r(v),&
		\text{for $u \ne v$.}
\end{align*}
The position of $\inter u'$ and $\inter v'$ forces $\inter{y_2}'$ to be placed between
$\inter{x_2}'$ and $\inter{y_2}'$ in every extending representation, which forces a $2$-FAT
obstruction.

\heading{$\boldsymbol{k}$-FS obstructions.} The class of \emph{forced side obstructions} is shown in
Fig.~\ref{fig:overview_of_obstructions}b.  To create $\widetilde H_k$ from $H_k$, we add a vertex
$u$ adjacent to $y_k$ and $z_k$. The partial representation contains three pre-drawn vertices $x_k$,
$y_k$ and $u$. We have $\ell(y_k) \le \ell(u) \le r(y_k) < r(u)$ and $\inter{x_k}'$ is on the
left of $\inter{y_k}'$. 

\heading{$\boldsymbol{k}$-EFS obstructions.}
The class of \emph{extended forced side obstructions} is similar to $k$-FS obstructions; see
Fig.~\ref{fig:overview_of_obstructions}c. To create $\widetilde H_k$ from $H_k$, we add $u$ adjacent
to $x_k$, $y_k$, and $z_k$, and $v$ adjacent to $y_k$ and $z_k$. The partial representation contains
four pre-drawn vertices $y_k$, $z_k$, $u$ and $v$ pre-drawn as follows:
$$\ell(u) < \ell(v) < \ell(y_k) \le r(y_k) < \ell(z_k) \le r(z_k) < r(u) \le r(v).$$

\heading{$\boldsymbol{k}$-FB obstructions.}
The class of \emph{forced betweenness obstructions} is similar to $k$-BI with $u=v$; see
Fig.~\ref{fig:overview_of_obstructions}d. To create $\widetilde H_k$ from $H_k$, we add $u$
adjacent to $y_k$ and $z_k$. The partial representation $\calR'_{\widetilde H_k}$ contains three
pre-drawn vertices $x_k$, $z_k$, and $u$.  We have $\ell(u) < \ell(z_k) \le r(u) \le r(z_k)$ and
$\inter{x_k}'$ is pre-drawn on the left of $\inter u'$.

\heading{$\boldsymbol{k}$-EFB obstructions.}
The class of \emph{extended forced betweenness obstructions} is similar to both $k$-BI and
$k$-FB; see Fig.~\ref{fig:overview_of_obstructions}e. To create $\widetilde H_k$ from $H_k$, we add $u$
adjacent to $x_k$, $y_k$ and $z_k$, and $v$ adjacent to $y_k$ and $z_k$. The partial representation
contains four pre-drawn vertices $x_k$, $z_k$, $u$, and $v$.  We have
$$\ell(u) < \ell(x_k) \le r(x_k) < \ell(v) < \ell(z_k) \le r(u) \le r(z_k) < r(v).$$

\heading{$\boldsymbol{k}$-FDS obstructions.}
The class of \emph{forced different sides obstructions} is shown in
Fig.~\ref{fig:overview_of_obstructions}f. To create $\widetilde H_k$ from $H_k$, we add $u$ adjacent
to $x_k$, $y_k$, and $z_k$, and $v$ adjacent to $y_k$ and $z_k$. The partial representation
contains three pre-drawn vertices $y_k$, $u$ and $v$ pre-drawn as follows:
$$\ell(u) < \ell(y_k) \le \ell(v) \le r(y_k) < r(u) \le r(v).$$

\heading{$\boldsymbol{k}$-EFDS obstructions.}
The class of \emph{extended forced different sides obstructions} is similar to $k$-FDS
obstructions; see Fig.~\ref{fig:overview_of_obstructions}g. To the construction of
$k$-FDS, we further add $w$ adjacent to $y_k$ and $z_k$.  The partial representation
contains four pre-drawn vertices $y_k$, $u$, $v$ and $w$ as follows:
$$\ell(u) < \ell(v) < \ell(y_k) \le \ell(w) \le r(y_k) < r(w) < r(u) \le r(v).$$

\heading{$\boldsymbol{k}$-FNS obstructions.} The class of \emph{forced nested side obstructions} is
constructed similarly as $k$-EFDS obstructions, but with $z_k$ pre-drawn instead of $y_k$; see
Fig.~\ref{fig:overview_of_obstructions}h. In $\calR'_{\widetilde H_k}$, we have
$$\ell(u) < \bigl\{\ell(v), \ell(w)\bigr\} \le \ell(z_k) \le r(w) < r(u) \le r(v),$$
where $\ell(v)$ and $\ell(w)$ are ordered arbitrarily.

\begin{figure}[t!]
\centering
\includegraphics{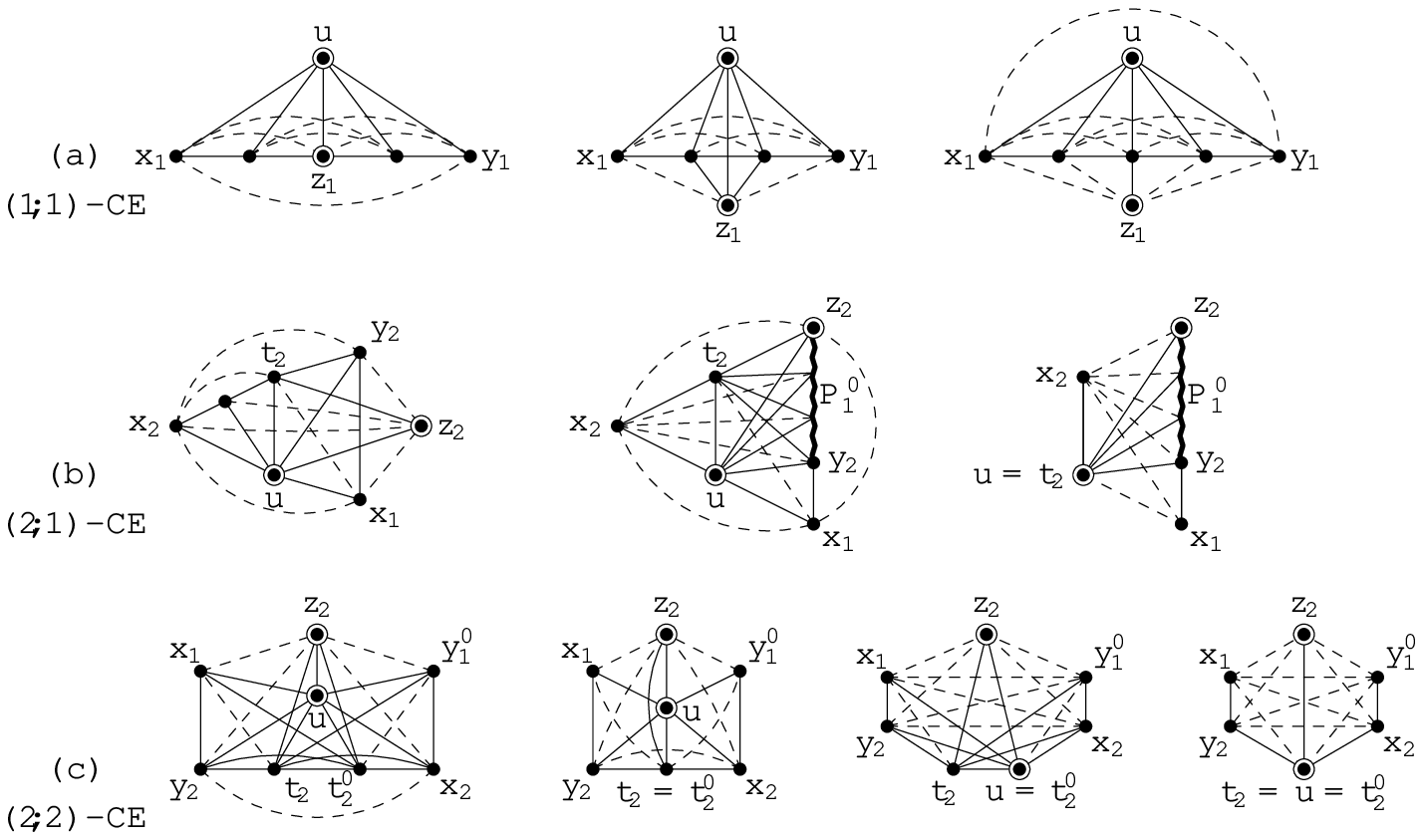}
\caption[]{(a) There are three minimal $(1,1)$-CE obstructions, and they are all finite graphs. The
reason is that if, say, a path $P_1$ was very long, we could replace $x_1$ by an inner vertex of
$P_1$; so such an obstruction would not be minimal.\\
(b) There are three minimal $(2,1)$-CE obstructions. Due to the minimality, $P_2$ is a path of length at most
two.  When $x_2$ and $t_2$ are non-adjacent, then $P'_1 = y_2t_2z_2$ is a path
avoiding $N[x_2]$. When $x_2$ and $t_2$ are adjacent, there are two cases, namely, $u \ne t_2$, and
$u = t_2$.  The path $P'_1$ is of length at least two.\\
(c) There are four minimal $(2,2)$-CE obstructions, and they are all finite graphs. Necessarily
$x_2t_2$ and $y_2t'_2$ are edges, and we can choose $x_2$ in such a way that it is adjacent to
$y'_1$ (similarly, $y_2$ is adjacent to $x_1$). There are four different graphs, because the
vertices $u$, $t_2$ and $t'_2$ might be distinct or not.}
\label{fig:kl_CE_obstructions}
\end{figure}

\heading{$\boldsymbol{(k,\ell)}$-CE obstructions.}
The class of \emph{covered endpoint obstructions} is created from a $k$-FAT obstruction glued
to an $\ell$-FAT obstruction; see Fig.~\ref{fig:overview_of_obstructions}i.  To create
$\widetilde H_{k,\ell}$, we glue $H_k$ with $H'_\ell$.  We put $z_k = z'_\ell$, $x_k = y'_\ell$ and
$y_k = x'_\ell$, and some other vertices of these obstructions may be also shared.  We add $u$
adjacent to $x_k$, $y_k$, and $z_k$. The partial representation contains two pre-drawn intervals
$\inter{z_k}'$ and $\inter{u}'$ such that $\ell(u) < \ell(z_k) \le r(u) \le r(z_k)$. We always
assume that $k \ge \ell$.

In $(k,\ell)$-CE Lemma~\ref{lem:kl_ce}, we show that the only $(k,\ell)$-CE obstructions which are
minimal are $(k,1)$-CE obstructions and $(2,2)$-CE obstructions; so either $\ell = 1$, or $k = \ell
= 2$. For $k \ge 3$, a minimal $(k,1)$-CE obstruction consists of the subgraph $H_k$ of a $k$-FAT
obstruction together with a vertex $u$, where $u$ is either adjacent to all vertices of $H_k$, or $u
= t_k$.  The remaining $(k,\ell)$-CE obstructions, where $2 \ge k \ge \ell$, are depicted in
Fig.~\ref{fig:kl_CE_obstructions}.

\subsection{Proofs of Non-extendibility and Minimality}

We sketch proofs that the defined obstructions are non-extendible and minimal. This implies the
first part of Theorem~\ref{thm:repext_obstructions}.  We establish the harder implication in
Sections~\ref{sec:strategy}, \ref{sec:leaves}, \ref{sec:p_nodes}, \ref{sec:q_nodes}, and
\ref{sec:main_result}. 

\begin{lemma} \label{lem:correctness_k_FAT}
Every $k$-FAT obstruction is non-extendible and minimal.
\end{lemma}

\begin{proof}
We prove the claim by induction. For $k=1$, non-extendibility and minimality are clear. For
$k>1$, assume that $\inter{x_k}'$ is on the left of $\inter{y_k}'$, and $\inter{y_k}'$ is on the left of $\inter{z_k}'$. In every representation of
$k$-FAT, $\inter{t_k}$ covers $[\ell(y_k),\ell(z_k)]$.  We know that
$\inter{x_{k-1}}$ cannot be on the left of $\inter{x_k}'$, since $H_{k-1}$ is connected and $x_k$ is
non-adjacent to all vertices of $H_{k-1}$.  Therefore $\inter{x_{k-1}}$ has to be placed on the
right of $\inter{z_k}'$. We get a $(k-1)$-FAT obstruction, which is non-extendible by
the induction hypothesis.

It remains to argue the minimality. If one of $\inter{x_k}'$, $\inter{y_k}'$ and $\inter{z_k}'$ is
made free, we can place them in such a way that $\inter{z_k}$ is between $\inter{x_k}$ and
$\inter{y_k}$. This makes the partial representation extendible: It works for $k=1$, and for $k>1$,
we can place $x_{k-1}$ on the right of $\inter{y_k}$, which makes the induced $H_{k-1}$ extendible.
If we remove one of the vertices or induced paths, the argument is similar.\qed
\end{proof}

\begin{lemma} \label{lem:correctness_all}
The following obstructions are non-extendible and minimal:
\begin{packed_itemize}
\item SE, $k$-FS, $k$-EFS, $k$-FB, $k$-EFB, $k$-FDS, $k$-EFDS, and $k$-FNS obstructions,
\item $k$-BI obstructions for $k \le 2$,
\item $(k,\ell)$-CE obstructions where either $\ell = 1$, or $k=\ell=2$.
\end{packed_itemize}
\end{lemma}

\begin{proof}
For SE obstructions, the proof is trivial. For the remaining classes (aside $k$-BI and
$(k,\ell)$-CE), we proceed as follows. Non-extendibility follows from the fact that, in all cases,
$\inter{y_k}$ is forced to be placed between $\inter{x_k}$ and $\inter{z_k}$. To show minimality, we
use the minimality of $k$-FAT obstructions. Then it is easy to show that freeing any added pre-drawn
interval or removing any added vertex results in the possibility of placing $\inter{z_k}$ between
$\inter{x_k}$ and $\inter{y_k}$.

Consider a $k$-BI where $k \le 2$.  Non-extendibility follows from the fact that $\inter{y_k}$ has
to be placed between $\inter{x_k}'$ and $\inter{z_k}'$, thus forcing the k-FAT obstruction,
which is non-extendible by Lemma~\ref{lem:correctness_k_FAT}.  By removing a vertex or an induced
path of $H_k$, it becomes extendible as argued in Lemma~\ref{lem:correctness_k_FAT}. By freeing
$\inter{u}'$ or $\inter{z_k}'$, we can place $\inter{y_k}$ on the right of $\inter{z_k}'$ which
makes the partial representation extendible. By freeing $\inter{v}'$ or $\inter{x_k}'$, we can place $\inter{y_k}$ on the
left of $\inter{x_k}'$, which also makes it extendible because $k \le 2$ and $x_2$ is adjacent to
$t_2$.

For $(k,\ell)$-CE, either $\inter{y_k}$ is between $\inter{x_k}$ and $\inter{z_k}'$ (non-extendible
due to the $k$-FAT obstruction), or $\inter{y_\ell}$ is between $\inter{x_\ell}$ and
$\inter{z_\ell}'$ (non-extendible due to the $\ell$-FAT obstruction). Minimality is also easy:
Removing or freeing $u$ allows to place $\inter{z_k}'$ between $\inter{x_k}$ and $\inter{y_k}$,
which is extendible. And removing anything from one of the FAT obstructions allows one of the
orderings of $\inter{x_k}$ and $\inter{y_k}$ to be extendible.\qed
\end{proof}

We note that the list of minimal obstructions is unique. Indeed, every minimal obstruction itself
corresponds to a valid input, which cannot be obstructed by a distinct obstruction due to the
minimality. Therefore, it is not possible to construct a smaller list of minimal obstructions, or to
argue that if the partial representation contains a particular obstruction, then it also contains an
additional one.

%% file: maximal_cliques.tex
\section{Maximal Cliques and MPQ-trees} \label{sec:maximal_cliques}

In this section, we review well-known properties of interval graphs. First, we describe their
characterization in terms of orderings of maximal cliques. Then, we introduce two data structures to
deal with these orderings, namely, PQ-trees and MPQ-trees. Finally, we prove some simple structural
results concerning MPQ-trees.

\heading{Consecutive Orderings.}
Fulkerson and Gross~\cite{maximal_cliques} proved the following fundamental characterization of interval
graphs:

\begin{figure}[b!]
\centering
\includegraphics{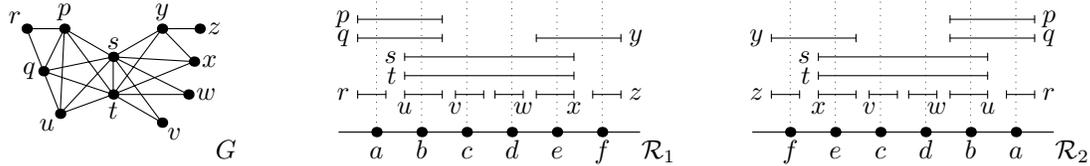}
\caption{An interval graph $G$ and two of its representations with different left-to-right orderings
$<$ of the maximal cliques. Some choices of clique-points are depicted on the real lines.}
\label{fig:fulkerson_gross}
\end{figure}

\begin{figure}[b!]
\centering
\includegraphics{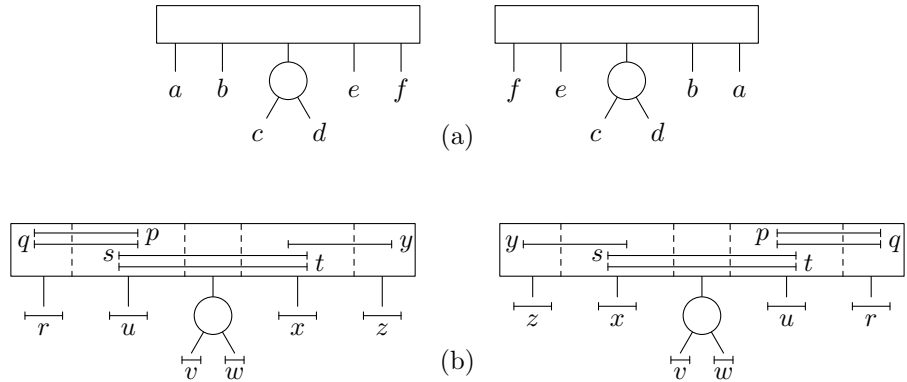}
\caption{(a) Two equivalent PQ-trees with frontiers $a < b < c < d < e < f$ and $f < e < c < d < b <
a$, respectively. In all figures, we denote P-nodes by circles and Q-nodes by rectangles. (b)~The
corresponding MPQ-trees with depicted sections.}
\label{fig:pq_trees}
\end{figure}

\begin{lemma}[Fulkerson and Gross~\cite{maximal_cliques}] \label{lem:maximal_cliques}
A graph is an interval graph if and only if there exists a linear ordering $<$ of its maximal
cliques such that, for each vertex, the maximal cliques containing this vertex appear consecutively.
\end{lemma}

We call an ordering of the maximal cliques satisfying the statement of
Lemma~\ref{lem:maximal_cliques} a \emph{consecutive ordering}. See Fig.~\ref{fig:fulkerson_gross}
for an example. The ordering $<$ from the statement is obtained by sweeping an interval
representation from left to right. By the Helly property, the intervals of every maximal clique have
a non-empty intersection. For all maximal cliques, these intersections are disjoint and ordered from
left to right. In the intersection of the intervals of a maximal clique $a$, we pick one point which
we call a \emph{clique-point} $\cp(a)$. The left-to-right ordering of these clique-points gives $<$.
On the other hand, given a consecutive ordering $<$, we place the clique-points from left to right
according to $<$ and construct an interval representation by placing each interval on top of its
clique-points and no others. This can be done because the ordering places the maximal cliques
containing each of the vertices consecutively.

\heading{PQ-trees.} Booth and Lueker~\cite{PQ_trees} designed a data structure called PQ-trees to
efficiently work with consecutive orderings of maximal cliques. A \emph{PQ-tree} $T$ is a rooted
tree. Its leaves are in one-to-one correspondence with the maximal cliques. Its inner nodes are of
two types: \emph{P-nodes} and \emph{Q-nodes}. Each P-node has at least two children, and each Q-node
has at least three. Further, for every inner node, a left-to-right ordering of its children is
given. Every PQ-tree $T$ represents one linear ordering $<_T$ of the maximal cliques called the
\emph{frontier} of $T$, which is the ordering of the leaves from left to right.

Every PQ-tree $T$ further represents other linear orderings. These orderings are frontiers of
equivalent PQ-trees. A PQ-tree $T'$ is \emph{equivalent} to $T$ if it can be constructed from $T$ by
a sequence of \emph{equivalent transformations} of two types: (i) an arbitrary reordering of the
children of a P-node, and (ii) a reversal of the order of the children of a Q-node.
Fig.~\ref{fig:pq_trees}a depicts two equivalent PQ-trees corresponding to the interval graph from
Fig.~\ref{fig:fulkerson_gross}.

Booth and Lueker proved that, for every interval graph, there exists a unique PQ-tree (up to
equivalence transformations) representing precisely the consecutive orderings of the maximal
cliques. In other words, this tree describes all possible interval representations of this interval
graph. 

A subtree of $T$ consists of a node and all its descendants. The subtrees of a node $N$ are those
subtrees having the children of $N$ as the roots. For a node $N$, let $T[N]$ denote the subtree of
$T$ with the root $N$.

\heading{MPQ-trees.}
For the purpose of this paper, we need more information about the way in which the vertices of the
interval graph are related to the structure of the PQ-tree. This additional information is contained
in the \emph{modified PQ-tree} (MPQ-tree), introduced by Korte and M\"ohring~\cite{korte_mohring}.
We note that the same idea is already present in the earlier paper of Colbourn and
Booth~\cite{colbourn_booth}.

The MPQ-tree is an augmentation of the PQ-tree in which the nodes of $T$ have assigned subsets of
$V(G)$ called \emph{sections}.  To a leaf representing a clique $a$, we assign one section $s(a)$.
Similarly, to each P-node $P$, we assign one section $s(P)$. For a Q-node $Q$ with $n$ children,
we have $n$ sections $s_1(Q),\dots,s_n(Q)$, each corresponding to one subtree of $Q$.

The section $s(a)$ has all vertices contained in the maximal clique $a$ and no other maximal clique.
The section $s(P)$ of a P-node $P$ has all vertices that are contained in all maximal cliques of
$T[P]$ and in no other maximal clique. Sections of Q-nodes are more complicated.
Let $Q$ be a Q-node with subtrees $T_1,\dots,T_n$. Let $x$ be a vertex contained only in maximal
cliques of $T[Q]$, and suppose that it is contained in maximal cliques of at
least two subtrees. Then $x$ is contained in every section $s_i(Q)$ such that some maximal clique of
$T_i$ contains $x$. Fig.~\ref{fig:pq_trees}b depicts the sections for the example in
Fig.~\ref{fig:fulkerson_gross}.

Korte and M\"ohring~\cite{korte_mohring} state the following properties:
\begin{packed_itemize}
\item Every vertex $x$ is placed in the sections of exactly one node of $T$. In the case of a
Q-node, it is placed in consecutive sections of this node.
\item For a Q-node $Q$, if $x$ is placed in a section $s_i(Q)$, then $x$ is contained in all cliques
of $T_i$.
\item Every section of a Q-node is non-empty. Moreover, two consecutive sections have a non-empty
intersection.
\item A maximal clique contains exactly those vertices contained in the sections encountered when we traverse the
tree from the corresponding leaf to the root.
\end{packed_itemize}

\heading{Structure of MPQ-trees.} Next, we show several structural properties used in
building minimal obstructions which are quite easy to prove:

\begin{lemma} \label{lem:different_sections}
Let $Q$ be a Q-node. Then $s_i(Q) \ne s_j(Q)$ for every $i \ne j$. Further, if $s_i(Q) \subsetneq
s_{i+1}(Q)$, then at least one section of $T_i$ is non-empty. 
\end{lemma}

\begin{proof}
If $s_i(Q) = s_j(Q)$, then we could exchange $T_i$ and $T_j$ and we would obtain a valid MPQ-tree
for $G$. Since $n \ge 3$, this yields a contradiction with the fact that the only possible
transformation of a $Q$-node is reverting the order of its children.

For the latter part, let $a$ and $b$ respectively be maximal cliques of a leaf in $T_i$ and a leaf in $T_{i+1}$. Then $a \setminus b \neq \emptyset$ and
every $x \in a \setminus b$ belongs to sections of $T_i$.\qed
\end{proof}

Let $N$ be a node of the MPQ-tree. By $G[N]$ we denote the subgraph induced by all the
vertices in the sections of the subtree rooted at $N$. Similarly, for a subtree $T'$, we denote the
subgraph induced by the vertices in its sections by $G[T']$.

\begin{lemma} \label{lem:connectivness}
Let $N$ be an inner node of an MPQ-tree.
\begin{packed_itemize}
\item[(i)] If $N$ is a Q-node, then $G[N]$ is connected.
\item[(ii)] If $N$ is a P-node, then $G[N]$ is connected if and only if $s(N)$ is non-empty.
Furthermore, for every child $T_i$ of $N$, the graph $G[T_i]$ is connected.
\end{packed_itemize}
\end{lemma}

\begin{proof}
(i)
It follows from the facts that the vertices in any section form a clique, and that
any two consecutive sections of $N$ have non-empty intersection.

(ii) The first statement is clear. For the second part, notice that if $G[T_i]$ was not connected, we
could permute the connected components of $G[T_i]$ arbitrarily with the other children of $N$.  Therefore
$T_i$ would not be a child of $N$, but $N$ would have one child per each connected component of
$G[T_i]$.\qed
\end{proof}

Let $Q$ be a Q-node and $i < j$. Let $x$ and $y$ be two vertices of $G[Q]$, where $x$ is either
in $T_i$, or $s_i(Q)$, and $y$ is either in $T_j$, or $s_j(Q)$. A path $P_{x,y}$ is called
\emph{$Q$-monotone} if all inner vertices of the path belong to the sections of $Q$, and their
leftmost/rightmost sections strictly increase.

\begin{lemma} \label{lem:monotone_paths}
Let $H$ be an induced subgraph of $G[Q]$ such that $x,y \in V(H)$ belong to one component.
Then every shortest path $P_{x,y}$ in $H$ is $Q$-monotone.
\end{lemma}

\begin{proof}
It is easy to see that any path from $x$ to $y$ that is not $Q$-monotone can be shortened.\qed
\end{proof}

Let $Q$ be a Q-node. Let $u$ be a vertex appearing in sections of $T[Q]$. If $u$ belongs to sections
of $Q$, let $s_u^\lft(Q)$ be the leftmost section of $Q$ containing $u$ and $s_u^\rt(Q)$ be the
rightmost one. If $u$ belongs to sections of a subtree $T_i$ of $Q$, we put $s_u^\lft(Q) =
s_u^\rt(Q) = s_i(Q)$.  If $s_u^\rt(Q)$ is on the left of $s_v^\lft(Q)$, then we say that \emph{$u$
is on the left of $v$} and \emph{$v$ is on the right of $u$}. Also, \emph{$u$ and $v$ are on the
same side of $w$} if they are both on the left of $w$, or both on the right of $w$. Similarly,
\emph{$v$ is between $u$ and $w$} if either $u$ is on the left and $w$ is on the right of $v$, or
$u$ is on the right and $w$ is on the left of $v$. For a maximal clique $a \in T_i$, we
say that \emph{$u$ is on the left of $a$} when $s_u^\rt(Q)$ is on the left of $s_i(Q)$, and
similarly the other relations.

\heading{Non-adjacencies of Maximal Cliques.} Maximal cliques of interval
graphs have the following special property, which we use to build minimal obstructions.

\begin{lemma} \label{lem:non_adjacent}
Let $H$ be a connected subgraph of an interval graph and let $c$ be a maximal clique with no vertex
in $V(H)$. There exists $x \in c$ non-adjacent to all vertices of $V(H)$.
\end{lemma}

\begin{proof}
Consider an interval representation $\calR$. It places all intervals of $H$ to one side of
$\cp(c)$, say on the left. Let $x$ be the interval of $c$ having the rightmost left endpoint. If $x$
is adjacent to some vertex $y \in V(H)$, then every vertex of $c$ is adjacent to $y$.  Since $c$ is
maximal, it follows that $y \in c$, contradicting the assumption. So $x$ is non-adjacent to all
vertices of $V(H)$.\qed
\end{proof}

%% file: interval_orders.tex
\section{Characterizing Extendible Partial Representations by Maximal Cliques} \label{sec:interval_orders}

In this section, we explain the characterization of extendible partial representations due to
Klav\'{\i}k et al.~\cite{kkosv}. 

\begin{figure}[t!]
\centering
\includegraphics{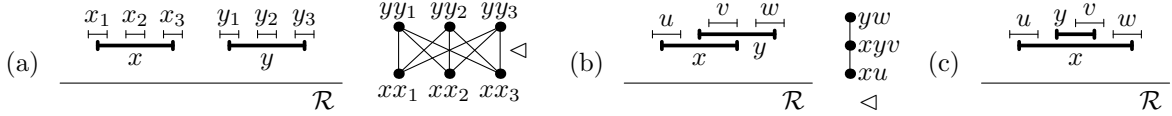}
\caption{Possible relative positions of pre-drawn intervals $\inter x'$ and $\inter y'$, and some
examples of the Hasse diagrams of the posed constraints. The consecutive ordering of the maximal
cliques of every extending representation has to extend $\wlt$. (a) All maximal cliques containing
$x$ have to be on the left of those containing $y$. (b) All maximal cliques containing $x$ have to
be on the left of those containing both $x$ and $y$, which are on the left of those containing only
$y$. (c) An inclusion of pre-drawn intervals poses no constraints. A maximal clique containing only
$x$ can be either on the left, or on the right of the maximal cliques containing both $x$ and $y$.}
\label{fig:simple_examples}
\end{figure}

\heading{Restricting Clique-points.} Suppose that there exists a representation $\calR$ extending
$\calR'$. Then $\calR$ gives some consecutive ordering $<$ of the maximal cliques from left to
right. We want to show that the pre-drawn intervals give constraints in the form of a partial
ordering $\wlt$.  Fig.~\ref{fig:simple_examples} illustrates these constraints given by a pair of
pre-drawn intervals. By generalizing it to all pre-drawn intervals, we get the partial ordering
$\wlt$.

\begin{figure}[b!]
\centering
\includegraphics{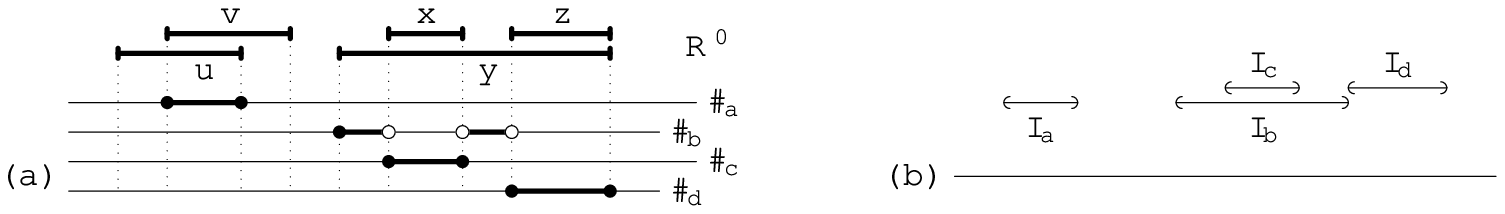}
\caption{(a) Four maximal cliques $a$, $b$, $c$, and $d$ with $P(a) = \{u,v\}$, $P(b) = \{y\}$, $P(c) =
\{x,y\}$, and $P(d) = \{y,z\}$. The possible positions $\da_a$, $\da_b$, $\da_c$, and $\da_d$ of
their clique-points are illustrated. (b) The corresponding open intervals $I_a$, $I_b$, $I_c$, and
$I_d$.}
\label{fig:restricting_cliques}
\end{figure}

For a maximal clique $a$, let $P(a)$ denote the set of all pre-drawn intervals that are contained in
$a$. Recall that a clique-point $\cp(a)$ is some point chosen from the intersection of all intervals
of $a$ in the representation $\calR$. Then $P(a)$ restricts the possible position of the
clique-point $\cp(a)$ to only those points $x$ of the real line which are covered in $\calR'$ by the
pre-drawn intervals of $P(a)$ and no others.  We denote the set of these admissible positions by
$\da_a$. Formally:
$$\da_a = \bigl\{ x : \text{$x \in \mathbb R$ and $x \in \inter u' \iff u \in P(a)$}\bigr\};$$
for examples see Fig.~\ref{fig:restricting_cliques}a. Equivalently, $\da_a$ is defined in~\cite{bko} as
\begin{equation} \label{eq:clique_point}
\da_a = \Bigl(\!\!\bigcap_{u \in P(a)}\!\!\! \inter u'\Bigr) \setminus
	\Bigl(\!\!\bigcup_{v \notin P(a)}\!\!\! \inter v'\Bigr).
\end{equation}

We are interested in the extremal points of $\da_a$. By $\clql(a)$ (resp.~$\clqr(a)$),
we denote the infimum (resp.~the supremum) of $\da_a$.  We use an open interval $I_a =
(\clql(a),\clqr(a))$ to represent $\da_a$. We note that this does not imply that $\da_a$ contains
all points between $\clql(a)$ and $\clqr(a)$; see $\da_b$ in Fig.~\ref{fig:restricting_cliques}.
Notice that when $P(a) = \emptyset$, then $I_a = \mathbb R$.

\heading{The Interval Order $\wlt$.} For two distinct maximal cliques $a$ and $b$, we write $a \wlt
b$ if $\clqr(a) \le \clql(b)$, or in other words, if $I_a$ is on the left of $I_b$. We put $a \wlt
a$ when $\da_a = \emptyset$. The definition of $\wlt$ is quite natural, since $a \wlt b$ implies
that every extending representation $\calR$ has to place $\cp(a)$ to the left of $\cp(b)$. For
instance, in Fig.~\ref{fig:restricting_cliques}, we get that $a \wlt b \wlt d$ and $a \wlt c \wlt
d$, but $b$ and $c$ are incomparable.

We note that the ordering $\wlt$ is a so-called \emph{interval order} represented by open intervals
$I_a$. The reason is that $a \wlt b$ if and only if $I_a$ and $I_b$ are disjoint and $I_a$ is on the
left of $I_b$. Interval orders are studied in the context of time constraints and have many
applications; for instance, see~\cite{fishburn}.

\begin{lemma}[Klav\'{\i}k et al.~\cite{kkosv}] \label{lem:repext_char}
A partial representation $\calR'$ is extendible if and only if there exists a consecutive ordering of
the maximal cliques that extends $\wlt$.
\end{lemma}

It is obvious that the constraints posed by $\wlt$ are necessary. The other implication is proved
in~\cite{kkosv} as follows. Suppose that $<$ is a consecutive ordering extending $\wlt$. We place
the clique-points greedily from left to right according $<$. We place each clique-point on the right
of the previously placed ones, as far to the left as possible.
It is shown that if this greedy procedure fails, then either $<$ does not extend $\wlt$, or it is
not a consecutive ordering.

\heading{Overlaps.}
In this paper, we show one additional crucial property of $\wlt$. We say that a pair of intervals
$I_a$ and $I_b$ \emph{single overlaps} if $I_a \ne I_b$ and either $\clql(a) \le \clql(b) < \clqr(a)
\le \clqr(b)$, or $\clql(b) \le \clql(a) < \clqr(b)
\le \clqr(a)$.

\begin{lemma} \label{lem:single_overlap}
No pair of intervals $I_a$ and $I_b$ single overlaps.
\end{lemma}

\begin{proof}
Assume without loss of generality that $\clql(a) \le \clql(b)$. If $\clqr(a) \le \clql(b)$, then
$I_a$ and $I_b$ are disjoint and do not single overlap. Suppose now that $\clql(a) \le \clql(b) <
\clqr(a)$. Since all intervals of $P(a)$ cover $[\clql(a),\clqr(a)]$, we get $P(a) \subseteq P(b)$.

The position of $\clqr(a)$ can be defined as a result of two distinct situations:
\begin{packed_itemize}
\item If some pre-drawn interval of $P(a)$ ends in $\clqr(a)$, then $\clqr(b) \le \clqr(a)$, since
the same pre-drawn interval is contained in $P(b)$. 
\item Otherwise, there exists a sequence of pre-drawn intervals not contained in $P(a)$ that covers
the whole portion between $\clqr(a)$ and the leftmost right endpoint of the intervals of $P(a)$.
The left endpoints of these intervals are on or to the right of $\clqr(a)$. Since the left endpoints
of the intervals in $P(b)$ are to the left of $\clqr(a)$, the pre-drawn intervals of the sequence
are not contained in $P(b)$. Thus, $\clqr(b) \le \clqr(a)$.
\end{packed_itemize}
In both cases, $\clql(a) \le \clql(b) \le \clqr(b) \le \clqr(a)$, so $I_b$ is contained in $I_a$.\qed
\end{proof}

If no single overlaps are allowed, every pair of intervals is either disjoint, or one interval is
contained in the other (possibly the intervals are equal). This type of interval orderings is very
simple and has not been much studied. We note that graphs having interval representations with no single
overlaps are called \emph{trivially perfect}. By examining the above proof, we get the following
useful result:

\begin{lemma} \label{lem:open_intervals}
If $I_a \subseteq I_b$, then $P(a) \supseteq P(b)$. Further, strict containments correspond
to strict inclusions.\qed
\end{lemma}

If $I_a$ and $I_b$ are disjoint, then we only know that at least one of the sets $P(a) \setminus
P(b)$ and $P(b) \setminus P(a)$ is non-empty. They both might be non-empty, or the sets $P(a)$ and
$P(b)$ might be in inclusion. See Fig.~\ref{fig:restricting_cliques} for examples.

\heading{Sliding Lemma.} We introduce some notation. We denote by $P^{\endlft}(a)$ and
$P^{\startrt}(a)$ respectively the subsets of $P(a)$ containing the pre-drawn intervals with
left-most right endpoints, and with right-most left endpoints. If $u \in P^{\endlft}(a)$ and $v \in
P^{\startrt}(a)$, then $\inter u' \cap \inter v' = \bigcap_{w \in P(a)} \inter w'$, thus $I_a$ is a
subinterval of $\inter u' \cap \inter v'$.

\begin{figure}[b!]
\centering
\includegraphics{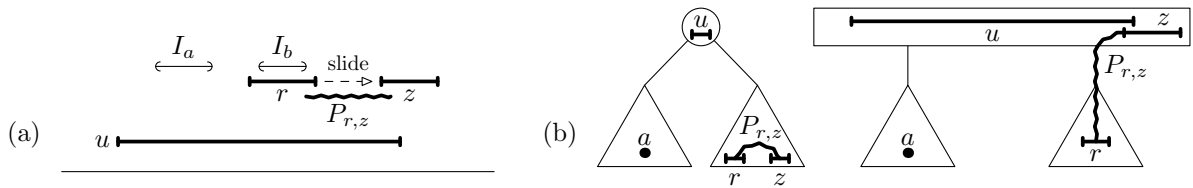}
\caption{(a) With the assumption satisfied, we can slide $r$ to $z$ which covers $r(u)$.
(b) The relative positions in the MPQ-trees of $z$ and $r$ with respect to $a$ are the same.}
\label{fig:sliding_lemma}
\end{figure}

Single overlaps of pre-drawn intervals pose more constraints than containment (see
Fig.~\ref{fig:simple_examples}b and c). Therefore, single overlaps are more powerful in building
obstructions. The following lemma states that, under some assumptions, we can turn a containment of
pre-drawn intervals into a single overlap of other pre-drawn intervals; see
Fig.~\ref{fig:sliding_lemma}.

\begin{lemma}[Sliding] \label{lem:sliding}
Let $I_a$ be on the left of $I_b$, $P(a) \subsetneq P(b)$ and $r \in P(b) \setminus P(a)$.
\begin{packed_enum}
\item[(i)]
There exists a pre-drawn interval $\inter z'$ on the right of $I_a$ covering $r(u)$, for $u \in
P^{\endlft}(a)$. Further, there exists an induced path $P_{r,z}$ from $r$ to $z$ whose vertices
are all pre-drawn and not contained in $P(a)$.
\item[(ii)]
Consider the smallest subtree having $a$ and the sections containing $r$. If the root of this subtree is
a P-node, then $z$ and $r$ are contained in the same subtree. If the root is a Q-node, then $z$ and $r$
appear on the same side of $a$.
\end{packed_enum}
\end{lemma}

\begin{proof}
(i) By the assumptions and the definition of $I_a$, we get that $(\clqr(a),r(u)]$ is not empty and
all points in $(\clqr(a),r(u)]$ are covered by pre-drawn intervals not in $P(a)$. Among these
intervals, we choose $z$ covering $r(u)$. Since $r$ is also one of these intervals, we can construct
an induced path from $r$ to $z$, consisting of pre-drawn intervals not in $P(a)$.

(ii) It follows from the existence of $P_{r,z}$ not contained in $a$.\qed
\end{proof}

We note that possibly $r = z$. The above lemma is repeatedly used for constructing minimal
obstructions. The general idea is the fact that $\inter r'$ properly contained inside $\inter u'$
restricts the partial representation less than $\inter z'$ covering $r(u)$. The lemma
says that we can assume that such $z$ exists and use it instead of $r$.

\heading{Flip Operation.} We say that we \emph{flip} the partial representation \emph{vertically}
when we map every $x \in \mathbb R$ to $-x$. This reverses the ordering $\wlt$.  Clearly, there
exists an obstruction in the original partial representation if and only if the flipped obstruction
is present in the flipped partial representation. The purpose of this operation is to decrease the
number of cases in the proofs.

%% file: strategy.tex
\section{Strategy for Finding Minimal Obstructions} \label{sec:strategy}

In this section, we describe the general strategy to show that every non-extendible partial
representation contains one of the obstructions described in
Section~\ref{sec:definition_of_obstructions}. 

For any two disjoint subtrees $T_i$ and $T_j$ of the MPQ-tree $T$, we write $T_i \wlt T_j$ if and
only if there exist cliques $a \in T_i$ and $b \in T_j$ such that $a \wlt b$. In this situation, the
maximal cliques of $T_i$ are forced to appear on the left of the maximal cliques of $T_j$.

\heading{Testing Extendibility by MPQ-tree Reordering.}
Recall that a MPQ-tree $T$ represents all feasible orderings of the maximal cliques of a given
interval graph $G$. By Lemma~\ref{lem:repext_char}, a partial representation is extendible if and
only if there exists a reordering $T'$ of $T$ such that the frontier of $T'$ extends $\wlt$. This
condition can be tested by the following algorithm (see~\cite{kkosv}).

\begin{figure}[b!]
\centering
\includegraphics{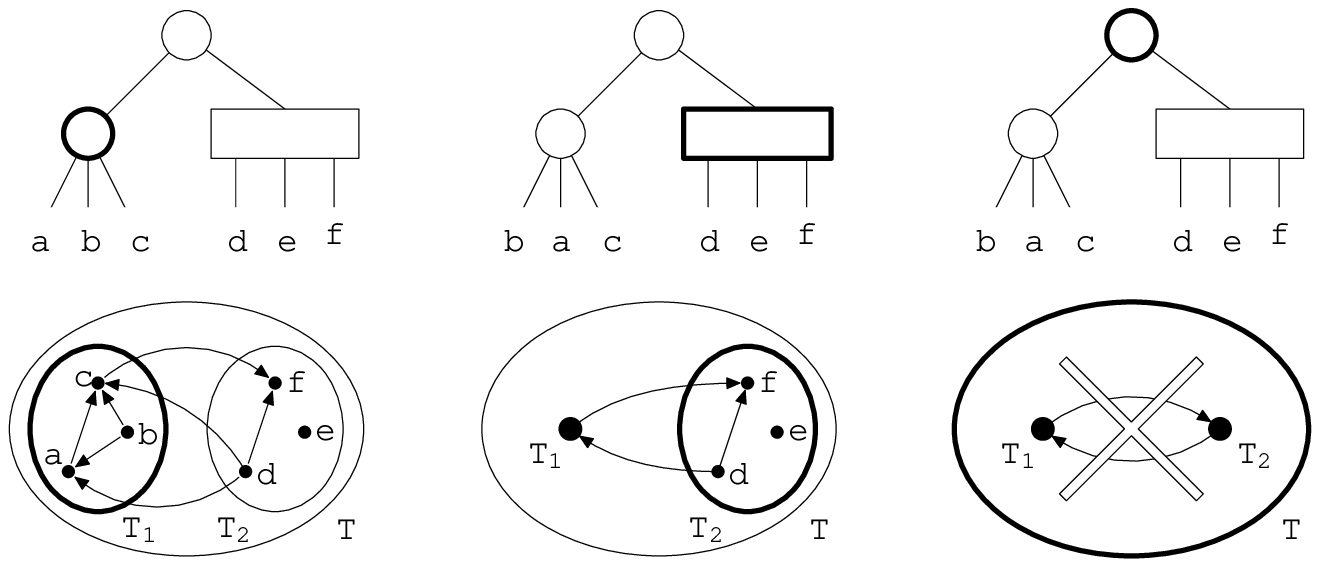}
\caption[]{This example is from~\cite{kkosv}, and it shows from left to right the way in which the reordering
algorithm works. We depict comparable pairs of maximal cliques by directed edges. The processed
trees are contracted into vertices.

\hskip 1em First, we reorder the highlighted P-node on the left. The subdigraph induced by $a$, $b$
and $c$ is ordered $b \rightarrow a \rightarrow c$. We contract this subtree $T_1$ into a vertex.
Next, we keep the order of the highlighted Q-node and contract its subtree $T_2$ into a vertex. When
we reorder the root P-node, the algorithm finds a two-cycle between $T_1$ and $T_2$, and outputs
``no''.}
\label{fig:reordering_example}
\end{figure}

We process the MPQ-tree $T$ from the bottom to the root. When a P-node is processed, we test whether
there exists a linear extension of $\wlt$ on its subtrees. It exists if and only if $\wlt$ induced
on the subtrees of the node is acyclic. Thus, if there is a cycle, the MPQ-tree cannot be reordered
according to $\wlt$. When a Q-node is processed, there are two possible orderings of its subtrees,
and we check whether any of them is compatible with $\wlt$. The partial representation is extendible
if and only if all nodes can be reordered in this manner. See Fig.~\ref{fig:reordering_example} for
an example. 

A node that cannot be reordered is called \emph{obstructed}. A set of maximal cliques
\emph{creates} an obstruction if the ordering of this set in $\wlt$ makes the node
obstructed.

\begin{figure}[t!]
\centering
\includegraphics{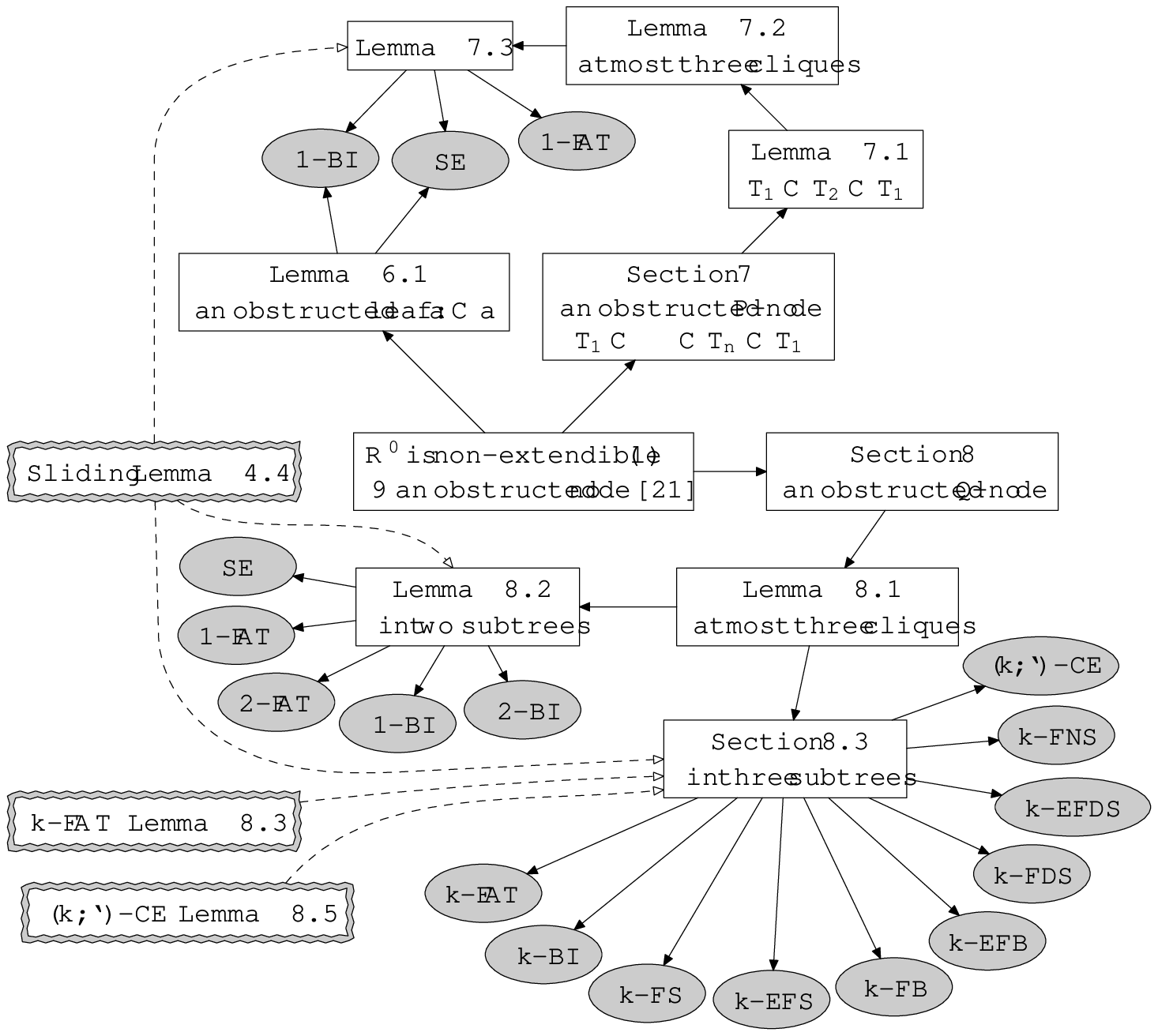}
\caption{Overview of the main steps of the proof, it starts in the middle. The obtained
obstructions are highlighted in gray, and three tools are depicted with highlighted borders. The
most involved case is in Section~\ref{sec:q_nodes_three_subtrees}.}
\label{fig:diagram}
\end{figure}

\heading{Strategy.} Suppose that a partial representation $\calR'$ is non-extendible.
From~\cite{kkosv}, we know that there exists an obstructed node in the MPQ-tree. We divide the
argument into three cases, according to the type of this node: \emph{an obstructed leaf}
(Section~\ref{sec:leaves}), \emph{an obstructed P-node} (Section~\ref{sec:p_nodes}), and \emph{an
obstructed Q-node} (Section~\ref{sec:q_nodes}).  Figure~\ref{fig:diagram} shows an overview of the
proof.

First, we argue that there exist at most three maximal cliques creating an obstruction.
Then, we consider their positions in the MPQ-tree and their open intervals from the definition of
$\wlt$. We use tools of Sections~\ref{sec:maximal_cliques} and~\ref{sec:interval_orders} to derive
positions of several pre-drawn intervals forming one of the obstructions.

In Section~\ref{sec:q_nodes_k_fat}, we prove a key tool called $k$-FAT Lemma~\ref{lem:k_fat}:
If three non-adjacent vertices $x_k$, $y_k$, and $z_k$ are pre-drawn in an order that is
different from their order in the sections of a Q-node, then they induce a $k$-FAT obstruction.
The proof is done by induction for $k$, and it explains why complicated obstructions are
needed.

%% file: leaves.tex
\section{Obstructed Leaves} \label{sec:leaves}

Suppose that some clique-point $a$ cannot be placed, so $\da_a = \emptyset$.  In terms of
$\wlt$, we get $a \wlt a$. Since $\wlt$ is a strict partial ordering, this already makes the partial
representation non-extendible.

\begin{lemma}[The leaf case] \label{lem:leaf}
If a leaf is obstructed, then $G$ and $\calR'$ contain an SE, or 1-BI obstruction.
\end{lemma}

\begin{proof}
We name the vertices as in the definition of the 1-BI obstructions. Suppose that the leaf corresponds
to a maximal clique $a$ such that $\da_a = \emptyset$.  

Let $u \in P^{\endlft}(a)$ and $v \in P^{\startrt}(a)$ (possibly $u = v$). Since $I_a$ is a
subinterval of $\bigcap_{w \in P(a)} \inter w'$ and $\da_a = \emptyset$, every point of
$[\ell(v),r(u)]$ is covered by some pre-drawn interval not contained in $P(a)$. Let $\inter{x_1}'$
be one such interval covering $\ell(v)$ and let $\inter{z_1}'$ be one such interval covering $r(u)$
(again, possibly $x_1=z_1$); see Fig.~\ref{fig:leaf}. Let $P_1$ be a shortest path from $x_1$ to
$z_1$ consisting of pre-drawn intervals not in $P(a)$.

\begin{figure}[t!]
\centering
\includegraphics{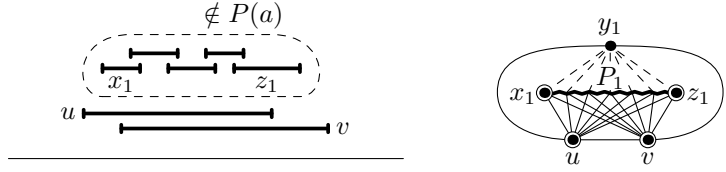}
\caption{A construction leading to a 1-BI obstruction.}
\label{fig:leaf}
\end{figure}

We prove that the relative pre-drawn position of $u$, $v$, $x_1$, and $z_1$ makes the partial
representation non-extendible. The maximal clique $a$ does not contain any vertex of $P_1$.  Since
the vertices of $P_1$ induce a connected subgraph, by Lemma~\ref{lem:non_adjacent} there exists $y_1
\in a$ which is non-adjacent to all vertices of $P_1$.  Hence, these (at most) five vertices
together with $P_1$ create either a 1-BI obstruction (when $\ell(u) < r(v)$), or an SE obstruction
(when $\ell(u)=r(v)$, for which $x = x_1 = z_1$ and we can free it). We note that this obstruction might not be minimal,
in which case we can remove some vertices and get one of the minimal obstructions illustrated
in Fig.~\ref{fig:SE} and~\ref{fig:BI_obstructions}.\qed
\end{proof}

%% file: p_nodes.tex
\section{Obstructed P-nodes} \label{sec:p_nodes}

If a P-node is obstructed, then it has some subtrees $T_1,\dots,T_n$ forming the cycle $T_1 \wlt T_2
\wlt \cdots \wlt T_n \wlt T_1$.  We start by showing that the specific structure of $\wlt$ forces
the existence of a two-cycle, so we can assume that $n=2$.

\begin{lemma} \label{lem:p_node_two_cycle}
If a P-node is obstructed, then it has two subtrees $T_1$ and $T_2$ such that $T_1 \wlt T_2 \wlt
T_1$.
\end{lemma}

\begin{figure}[b!]
\centering
\includegraphics{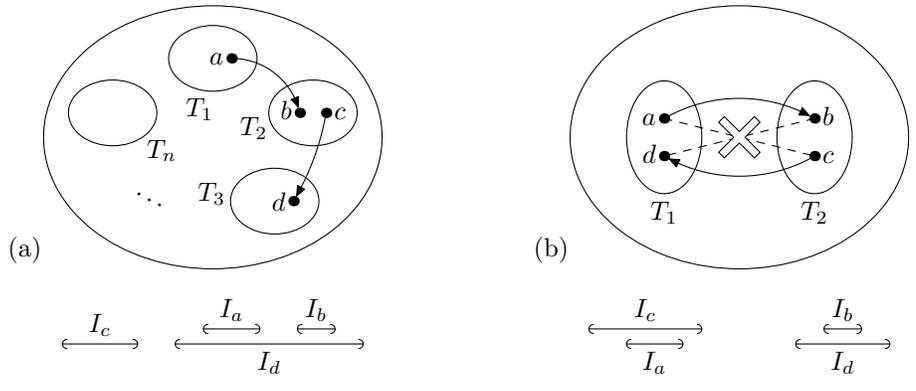}
\caption{(a) At the top, a shortest $n$-cycle of $\wlt$ on the children of a P-node. At the bottom, the
derived positions of the open intervals. (b) At the top, the four cliques involved in a two-cycle in $\wlt$.
The cliques $a$ and $c$ are incomparable, and so are $b$ and $d$. At the bottom, one of the four
possible configurations of the open intervals.}
\label{fig:P_node_cycles}
\end{figure}

\begin{proof}
The proof is illustrated in Fig.~\ref{fig:P_node_cycles}a.  Let $T_1 \wlt \cdots \wlt T_n \wlt T_1$
be a shortest cycle for the obstructed P-node. To get a contradiction, we assume $n \ge 3$. Since
$T_1 \wlt T_2$, there exist $a \in T_1$ and $b \in T_2$ such that $a \wlt b$.  Similarly, there
exist $c \in T_2$ and $d \in T_3$ such that $c \wlt d$.  We know that $I_a$ is on the left of $I_b$,
and $I_c$ is on the left of $I_d$. We analyze the remaining relative positions.

First, $I_d$ is not on the right of $I_a$, since otherwise $T_1 \wlt T_3$ and a shorter cycle would
exist. Additionally, $I_d$ is not on the left of $I_b$, since we would get $T_3 \wlt T_2$, and $T_2$
and $T_3$ would form a two-cycle. According to Lemma~\ref{lem:single_overlap}, no single overlaps of
open intervals are allowed, so $I_d$ necessarily contains both $I_a$ and $I_b$; see
Fig.~\ref{fig:P_node_cycles}a. Therefore, $I_c$ is on the left of $I_a$, so $T_2 \wlt T_1$ and we
get a two-cycle.\qed
\end{proof}

To create a two-cycle, at most four cliques are enough.  Aside from Lemma~\ref{lem:single_overlap},
so far we have not used that $\wlt$ arises from a partial interval representation. Next, we use
properties of the MPQ-tree.

\begin{lemma} \label{lem:p_node_three_cliques}
A two-cycle $T_1 \wlt T_2 \wlt T_1$ is created by at most three cliques.
\end{lemma}

\begin{proof}
The proof is depicted in Fig.~\ref{fig:P_node_cycles}b. Suppose that this two-cycle is given by four
cliques $a,d \in T_1$ and $b,c \in T_2$ such that $a \wlt b$ and $c \wlt d$. Assume for
contradiction that no three of these cliques define the two-cycle, i.e., $a$ and $c$ are
incomparable, and so are $b$ and $d$. According to Lemma~\ref{lem:single_overlap}, $I_a \subseteq
I_c$ or $I_a \supseteq I_c$, and analogously for $I_b$ and $I_d$. In all of the four cases, $I_c$ is
on the left of $I_b$, and $I_d$ is on the right of $I_a$.

We look at the case where $I_a \subseteq I_c$ and $I_b \subseteq I_d$, as in
Fig.~\ref{fig:P_node_cycles}b. By Lemma~\ref{lem:open_intervals}, we have $P(c) \subseteq P(a)$.
Therefore, $P(c)$ contains no vertices from the sections of $T_2$. Similarly, $P(d) \subseteq P(b)$,
and $P(d)$ contains no vertices from the sections of $T_1$. Therefore $P(c) = P(d)$, which implies
$I_c = I_d$, a contradiction. The other cases can be analyzed similarly.\qed
\end{proof}

It remains to put these results together and characterize the possible obstructions.

\begin{lemma}[The P-node case] \label{lem:p_node}
If a P-node is obstructed, then $G$ and $\calR'$ contain an SE, $1$-FAT, or $1$-BI obstruction.
\end{lemma}

\begin{proof}
According to Lemma~\ref{lem:p_node_two_cycle}, the obstructed P-node has a two-cycle in $\wlt$.  By
Lemma~\ref{lem:p_node_three_cliques}, there are at most three maximal cliques defining this cycle.  First
assume that this cycle is defined by two cliques $a \in T_1, b \in T_2$ such that $a \wlt b \wlt a$.
According to the definition of $\wlt$, this implies that $I_a = I_b$, both of lenght zero.
Therefore $P(a) = P(b)$. Let $u \in P^{\endlft}(a)$ and $v \in P^{\startrt}(a)$ (possibly $u=v$); we
have that $\inter u' \cap \inter v'$ is a singleton.  Since $a$ and $b$ are two maximal cliques,
there exists $x \in a \setminus b$ and $y \in b \setminus a$. We get an SE obstruction.

It remains to deal with the case where three cliques define the two-cycle. Let $a,c \in T_1$ and $b
\in T_2$ such that $a \wlt b \wlt c$. We have three non-intersecting intervals whose left-to-right
order is $I_a$, $I_b$ and $I_c$. Since $I_a$ and $I_c$ are disjoint, one of the sets $P(a) \setminus
P(c)$ and $P(c) \setminus P(a)$ is non-empty. Without loss of generality, we assume that $P(a)
\setminus P(c)\neq \emptyset$. Let $p \in P(a) \setminus P(c)$; then $p$ belongs to sections of
$T_1$, and as a consequence $p \notin P(b)$. Therefore $\inter p'$ is on the left of $I_b$.  We
distinguish two cases.

\begin{figure}[b!]
\centering
\includegraphics{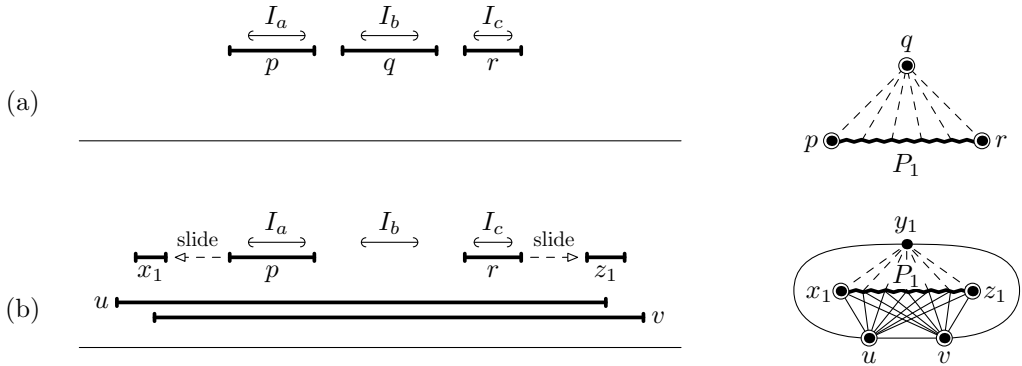}
\caption{The two cases of the proof of Lemma~\ref{lem:p_node}. (a) Case 1 leads to a 1-FAT
obstruction. (b) Case 2 leads to a 1-BI obstruction.}
\label{fig:P_node_two_cases}
\end{figure}

\emph{Case 1:} $P(b) \setminus P(c) \ne \emptyset$. We choose $q \in P(b) \setminus P(c)$. Then
$q$ belongs to sections of $T_2$, and $\inter q'$ is between $\inter p'$ and $I_c$. In
the next paragraph we show that there also exists $r \in P(c) \setminus P(b)$. Then $r$ belongs to sections of
$T_1$, and $\inter r'$ is on the right of $\inter q'$, as in Fig.~\ref{fig:P_node_two_cases}a. By
Lemma~\ref{lem:connectivness}(ii), $G[T_1]$ is connected; let $P_1$ be a shortest
path from $p$ to $r$ in $G[T_1]$. We obtain a $1$-FAT obstruction for $x_1 = p$, $y_1 = q$ and $z_1 = r$.

It remains to show that such $r$ exists. Suppose for contradiction that $P(c) \setminus P(b)
= \emptyset$. Since $P(c) \subsetneq P(b)$, no vertex of $P(c)$ appears in sections of $T_1$,
and we get $P(c) \subsetneq P(a)$. Consequently, every pre-drawn interval of $P(c)$ contains
$[\clql(a),\clqr(c)]$. The position of $I_c$ implies that every point of $[\clql(a),\clql(c))$ is
covered by some pre-drawn interval not contained in $P(c)$. In particular, there exists a path from
$p$ to $q$ consisting of such intervals. Since $p$ belongs to
sections of $T_1$ and $q$ belongs to sections of $T_2$, every path from $p$ to $q$ 
contains a vertex of the section of the P-node, or of a section above it; hence, the path contains a
vertex belonging to $c$. We obtain a contradiction.

\emph{Case 2:} $P(b) \setminus P(c) = \emptyset$. Then there exists $r \in P(c) \setminus P(b)$. 
We again observe that $\inter r'$ is on the right of $I_b$, as depicted in
Fig.~\ref{fig:P_node_two_cases}b. Furthermore, $P(b) \subseteq P(a) \cap P(c)$, so every
pre-drawn interval of $P(b)$ contains $[\clql(a),\clqr(c)]$.

We construct a 1-BI obstruction and we name the vertices as in the definition. Let $u \in
P^\endlft(b)$ and $v \in P^\startrt(b)$ (possibly $u = v$). Since $p$ does not necessarily cover
$\ell(v)$ and $r$ does not necessarily cover $r(u)$, we might not be able to construct a $1$-BI
obstruction with $x_1=p$ and $z_1=q$. We instead use Sliding Lemma~\ref{lem:sliding}. By applying it
(flipped) to $I_b$, $I_a$ and $p$, we obtain a pre-drawn interval $x_1$ covering $\ell(v)$ (possibly
$x_1 = p$). By applying it to $I_b$, $I_c$ and $r$, we obtain a pre-drawn interval $z_1$ covering
$r(u)$ (possibly $z_1 = r$). Furthermore, $x_1$ and $z_1$ belong to sections of $T_1$. Since
$G[T_1]$ is connected by Lemma~\ref{lem:connectivness}(ii), there exists a shortest path $P_1$ from
$x_1$ to $z_1$ containing no vertex of $b$. By Lemma~\ref{lem:non_adjacent}, there exists $y_1 \in
b$ non-adjacent to all vertices of $P_1$. We obtain a $1$-BI obstruction.\qed
\end{proof}

%% file: q_nodes.tex
\section{Obstructed Q-nodes} \label{sec:q_nodes}

Suppose that a Q-node with subtrees $T_1,\dots,T_n$ is obstructed. Then the two possible
orderings of this Q-node are not compatible with $\wlt$. Notice that at most four cliques are
sufficient to create the obstruction. We next prove that at
most three cliques are already sufficient.

\begin{lemma} \label{lem:q_node_three_cliques}
If a Q-node is obstructed, there exists an obstruction created by at most three maximal cliques.
\end{lemma}

\begin{proof}
Suppose that an obstruction is created by four cliques $a \in T_\alpha$, $b \in T_\beta$, $c \in
T_\gamma$ and $d \in T_\delta$ such that $\alpha < \beta$, $\gamma < \delta$, $a \wlt b$, and $c
\wgt d$. We know that $I_a$ is on the left of $I_b$, and $I_c$ is on the right of $I_d$. 
Notice that the four subtrees $T_\alpha$, $T_\beta$, $T_\gamma$ and $T_\delta$ are not necessarily
distinct. We classify all possible orderings $<$ of $\alpha$, $\beta$, $\gamma$, $\delta$ in two
general cases, namely, $\alpha \neq \gamma$ and $\alpha = \gamma$. In the first case, we may assume
without loss of generality that $\alpha < \gamma$.

\emph{Case 1: $\alpha < \gamma < \delta$} (see Fig.~\ref{fig:Q_node_three_cliques}a). Consider the
relative positions of $I_c$ and $I_d$ with respect to $I_a$. If $I_d$ is to the left of $I_a$, we
have $d \wlt a \wlt b$, and these three cliques already create an obstruction. If $I_c$ is to the
right of $I_a$, then we get $a \wlt c$ and $c \wgt d$, creating an obstruction. If neither happens,
then $I_c$ and $I_d$ are subintervals of $I_a$. Thus $c,d \wlt b$. If $\beta \le \gamma$, we have $a
\wlt b$ and $b \wgt d$, creating an obstruction. If $\beta > \gamma$, then $d \wlt c \wlt b$, which
also creates an obstruction.

\begin{figure}[b!]
\centering
\includegraphics{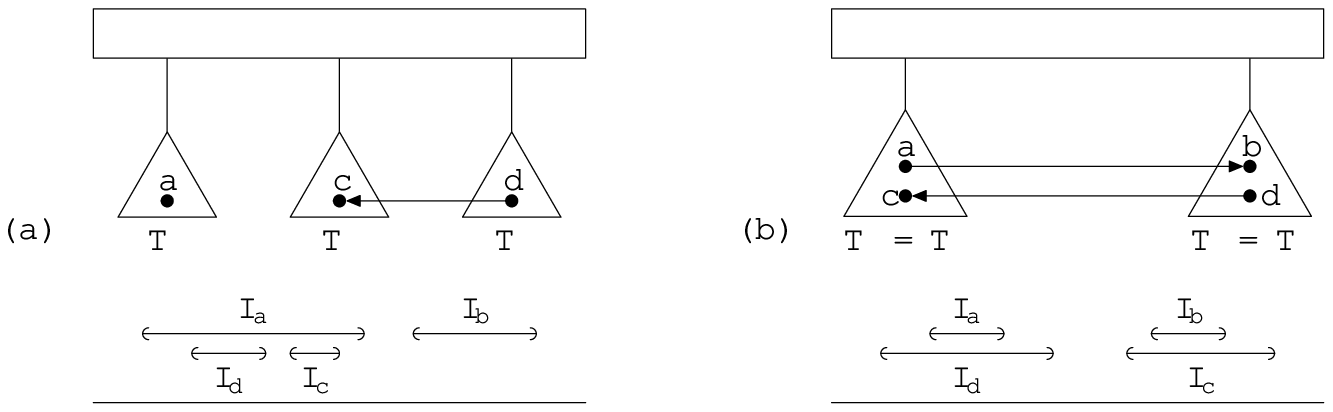}
\caption{Two cases of the proof of Lemma~\ref{lem:q_node_three_cliques}. The Q-node is depicted in
the top, while in the bottom we have the relative positions of the intervals.}
\label{fig:Q_node_three_cliques}
\end{figure}

\emph{Case 2: $\alpha = \gamma$} (see Fig.~\ref{fig:Q_node_three_cliques}b). If $I_c$ does not
intersect $I_b$, or $I_d$ does not intersect $I_a$, it is easy to see that three of the cliques
already create an obstruction. Suppose next that these intersections occur. Then $d \wlt b$. If
$\delta < \beta$ or $\beta < \delta$, it is again easy to show that three cliques are enough to
create an obstruction. It only remains to consider the case where $\alpha = \gamma < \beta =
\delta$.

Since the intervals $I_c$ and $I_a$ are non-intersecting, we may assume without loss of generality
that there exists $x \in P(a) \setminus P(c)$. This vertex $x$ belongs to sections of $T_\alpha$.
Thus $x \notin P(d)$, and we get that $I_a \subsetneq I_d$. By Lemma~\ref{lem:open_intervals}, $P(d)
\subsetneq P(a)$; in particular, $P(d)$ contains no private pre-drawn interval from sections of
$T_\beta$, and all pre-drawn intervals of $s_\beta(Q)$ are also contained in $s_\alpha(Q)$.

Since $P(d) \setminus P(b) = \emptyset$, there exists $y \in P(b) \setminus P(d)$ which is contained
in sections of $T_\beta$. We next apply the argument in the previous paragraph, and obtain $y \notin
P(c)$, $I_b \subsetneq I_c$, and $P(c) \subsetneq P(b)$. Consequently, $P(c)$ contains no private
pre-drawn intervals from sections of $T_\alpha$, and all pre-drawn intervals of $s_\alpha(Q)$ are
contained in $s_\beta(Q)$. We conclude that $P(c) = P(d)$ and $I_c = I_d$, which gives a
contradiction.\qed
\end{proof}

In summary, we can assume that a minimal obstruction involves at most three maximal cliques.  These
three cliques belong to either two or three different subtrees.

In the rest of the section, many figures describe positions of derived pre-drawn intervals in
sections of the Q-node and its subtrees; for instance Fig.~\ref{fig:Q_node_two_subtrees}. Some of
these intervals necessarily belong to sections of the Q-node, since they belong to maximal cliques
of several subtrees; for instance $t_2$ in Fig.~\ref{fig:Q_node_two_subtrees}. But for the remaining
intervals, it is not important to distinguish whether they belong to sections of the Q-node or one
of its subtrees, only their relative positions in the Q-node matter; for instance $q$ and $x_1$ in
Fig.~\ref{fig:Q_node_two_subtrees}.

\subsection{Cliques in Two Different Subtrees} \label{sec:q_nodes_two_subtrees}

In this section, we deal with the case where the maximal cliques belong to two different subtrees.

\begin{lemma}[The Q-node case, Two Subtrees] \label{lem:q_node_two}
If at most three cliques creating the obstruction belong to two different subtrees, then $G$ and
$\calR'$ contain an SE, $1$-FAT, $2$-FAT, $1$-BI, or $2$-BI obstruction.
\end{lemma}

\begin{proof}
The proof is similar to that of Lemma~\ref{lem:p_node}.  If two maximal cliques create an
obstruction, we can argue as in the first paragraph of the proof of Lemma~\ref{lem:p_node}, and we
obtain an SE obstruction.  It remains to deal with the case of three maximal cliques $a$, $b$, and
$c$.

We can assume that $a \wlt b \wlt c$ and that, for some $i < j$, we have $a,c \in T_i$ and $b \in
T_j$. Furthermore, without loss of generality, there exist $p \in P(a) \setminus P(c)$. Since $p$
belongs to sections of $T_i$, then $p \notin P(b)$, and thus $\inter p'$ lies to the left of $I_b$.
We distinguish two cases. 

\emph{Case 1:} $P(b) \setminus P(c) \ne \emptyset$. Then there exists $q \in P(b) \setminus P(c)$
such that $\inter q'$ lies between $\inter p'$ and $I_c$. Since $q$ is non-adjacent to $p$, it
belongs to sections of either $Q$ or $T_j$. Notice that in any case $s_q^\lft(Q)$ is on the right of
$s_i(Q)$.  Arguing as in Case 1 of the proof of Lemma~\ref{lem:p_node}, we observe that there exists
$r \in P(c) \setminus P(b)$. Furthermore, it follows that $\inter r'$ lies to the right of $\inter
q'$;  see Fig.~\ref{fig:Q_node_two_subtrees}a on the left.

\begin{figure}[b!]
\centering
\includegraphics{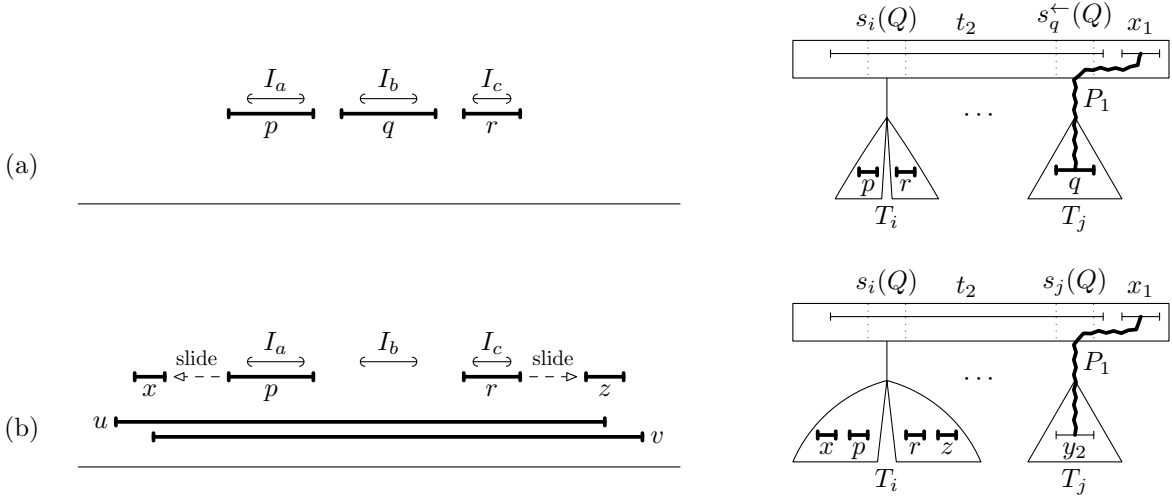}
\caption{(a) Case 1: The pre-drawn intervals and the situation in the
MPQ-tree for $s_i(Q) \subsetneq s_q^\lft(Q)$. (b) Case 2: The pre-drawn intervals and the
situation when there exists no path from $x$ to $z$ avoiding the vertices of $b$.}
\label{fig:Q_node_two_subtrees}
\end{figure}

If there exists a path $P_1$ from $p$ to $r$ avoiding $N[q]$, we get a $1$-FAT obstruction for $x_1
= p$, $y_1 = q$, $z_1 = r$ and $P_1$. By Lemma~\ref{lem:different_sections}, we know that $s_i(Q)
\ne s_q^\lft(Q)$. If $s_i(Q) \not\subseteq s_q^\lft(Q)$, then there exists some $w \in s_i(Q)
\setminus s_q^\lft(Q)$. Therefore, $P_1 = pwr$ is such a path. It remains to deal with the case
where no such path $P_1$ exists, which implies that $s_i(Q) \subsetneq s_q^\lft(Q)$; see
Fig.~\ref{fig:Q_node_two_subtrees}a on the right.

Consider the set $W = s_i(Q)$. Let $t_2$ be a vertex of $W$ whose section $s_{t_2}^\rt(Q)$ is
leftmost. Let $C$ be the component of $G[Q] \setminus W$ containing $q$. Since $s_q^\lft(Q)
\setminus W$ is non-empty, $C$ consists of the vertices of at least two subtrees of the Q-node. If
$t_2$ was adjacent to all vertices of $C$, it would be possible to flip the ordering of this
component, contradicting the fact that there are only two possible orderings for $Q$. Therefore,
$t_2$ is not adjacent to all vertices of $C$.  We choose $x_1 \in C \setminus N[t_2]$ whose section
$s_{x_1}^\lft(Q)$ is leftmost.  Let $P_1$ be a shortest path from $q$ to $x_1$ whose inner vertices
are adjacent to $t_2$. It follows that $x_2 = p$, $y_2 = q$, $z_2 = r$, $P_2 = x_2t_2$, $t_2$,
$x_1$, and $P_1$ define a $2$-FAT obstruction. (By Lemma~\ref{lem:monotone_paths}, all inner
vertices of $P_1$ are adjacent to $t_2$.)\qed

\emph{Case 2:} $P(b) \setminus P(c) = \emptyset$. Then there exists $r \in P(c) \setminus P(b)$.
Since $\inter r'$ lies on the right of $I_b$, the vertex $r$ is not contained in $a$ and it belongs
to sections of $T_i$. We use the same approach as in Case 2 of the proof of Lemma~\ref{lem:p_node}.
Since $P(b) \subseteq P(a) \cap P(c)$, every pre-drawn interval of $P(b)$ covers
$[\clql(a),\clqr(c)]$.  Let $u \in P^\endlft(b)$ and $v \in P^\startrt(b)$ (possibly $u=v$).

By applying Sliding Lemma~\ref{lem:sliding} twice, we get $x,z \notin P(b)$ such that $\inter x'$
covers $\ell(v)$ and $\inter z'$ covers $r(u)$; see Fig.~\ref{fig:Q_node_two_subtrees}b on the left.
Suppose that there exists a path $P_{x,z}$ from $x$ to $z$ avoiding all vertices of $b$. Let $x_1 =
x$, $z_1 = z$, and $P_1$ be a shortest path from $x_1$ to $z_1$ in $G[Q] \setminus b$.
By Lemma~\ref{lem:non_adjacent}, there exists $y_1 \in b$ non-adjacent to $P_1$. We
obtain a $1$-BI obstruction.

Suppose next that there is no path $P_{x,z}$ avoiding $b$.  We know that $x$ and $y$ belong to
sections of $T_i$, since there exist paths $P_{x,p}$ and $P_{r,z}$ avoiding $b$, from the above
applications of Sliding Lemma~\ref{lem:sliding}.  Since no path $P_{x,z}$ avoiding $b$ exists, we
have $s_i(Q) \subsetneq s_j(Q)$.  As in Case 1, let $W = s_i(Q)$, and let $t_2$ be a vertex of $W$
whose section $s_{t_2}^\rt(Q)$ is leftmost (possibly $t_2 = u$ or $t_2 = v$).  We again infer that
$t_2$ is not adjacent to all vertices of $C$, where $C$ is the component of $G[Q] \setminus W$
containing $b \setminus W$. We choose $x_1 \in C \setminus N[t_2]$ whose section $s_{x_1}^\lft(Q)$
is leftmost. Since $s_i(Q) \subsetneq s_j(Q) \subseteq b$, there exists $y_2 \in b$ non-adjacent to
$x$ and $z$. We get a 2-BI obstruction for $x_2 = x$, $y_2$, $z_2 = z$, $u$, $v$, a shortest path
$P_1$ from $y_2$ to $x_1$ in $C$, and $P_2 = x_2t_2$. (By Lemma~\ref{lem:monotone_paths}, all inner
vertices of $P_1$ are adjacent to $t_2$.)\qed
\end{proof}

\subsection{$k$-FAT and $(k,\ell)$-CE Lemmas} \label{sec:q_nodes_k_fat}

In this section, we give two tools for the case, analyzed in
Section~\ref{sec:q_nodes_three_subtrees}, where the three maximal cliques creating the obstruction belong to
three different subtrees. These tools give insight into the structure of the Q-nodes, and explain
the way in which complex obstructions such as $k$-FAT and $(k,\ell)$-CE obstructions are formed.

\heading{$\boldsymbol k$-FAT Lemma.} First, we present a useful lemma that allows to locate $k$-FAT
obstructions. The key idea of the proof is similar to Case 1 of the proof of
Lemma~\ref{lem:q_node_two}, but applied inductively for $k$. 

\begin{lemma}[$\boldsymbol k$-FAT] \label{lem:k_fat}
Let $Q$ be a Q-node with children $T_1,\dots,T_n$, and let $a$, $b$ and $c$ be three cliques of
$T[Q]$ contained respectively in $T_\alpha$, $T_\beta$ and $T_\gamma$, for $\alpha < \beta < \gamma$.
Let $x_k \in P(a)$, $y_k \in P(c)$ and $z_k \in P(b)$ be three disjoint pre-drawn intervals such
that $\inter{y_k}'$ is between $\inter{x_k}'$ and $\inter{z_k}'$. Then $G[Q]$ and
$\calR'[\{x_k,y_k,z_k\}]$ contain a $k$-FAT obstruction.
\end{lemma}

\begin{proof}
The proof, illustrated in Fig.~\ref{fig:q_nodes_k_fat_statement}, is by induction. We always denote
the vertices as in the definition of $k$-FAT obstructions. If we find a $1$-FAT or $2$-FAT
obstruction, the statement is true.  Otherwise, we recurse on a smaller part of the Q-node, where we
find a structure identical to a $(k-1)$-FAT obstruction, except for the fact that the vertex
$x_{k-1}$ is free. Together with some vertices in the remainder of the Q-node, we obtain a $k$-FAT
obstruction. We next provide the details.

\begin{figure}[b!]
\centering
\includegraphics{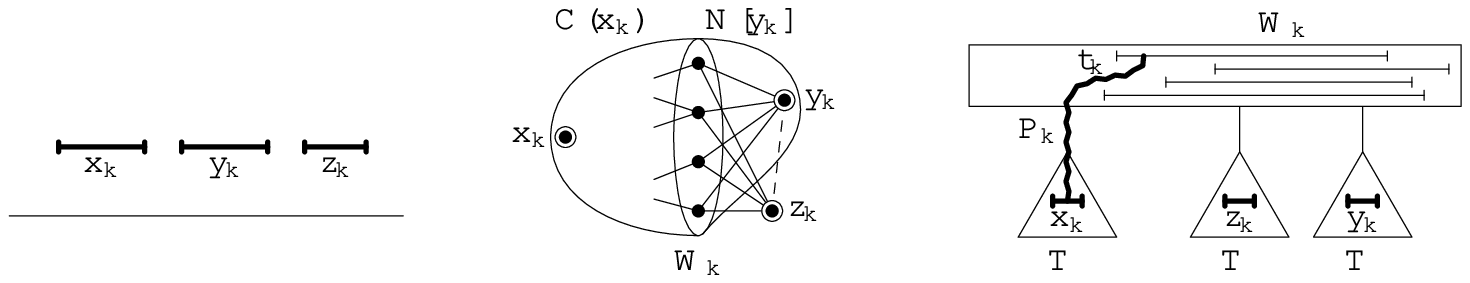}
\caption{On the left, the position of the pre-drawn intervals.  In
the middle, the construction of $W_k \subsetneq N[y_k]$ in $G[Q]$. On the right, the Q-node with the
three subtrees and the intervals of $W_k$ depicted in its sections.}
\label{fig:q_nodes_k_fat_statement}
\end{figure}

Let $k$ be some yet unspecified integer, determined by the recursion. We want to argue that $G[Q]$
contains a $k$-FAT obstruction because the ordering of $\inter{x_k}'$, $\inter{y_k}'$ and
$\inter{z_k}'$ is incorrect (in every representation, $\inter{z_k}$ is between $\inter{x_k}$ and
$\inter{y_k}$). Suppose that there exists a path from $x_k$ to $z_k$ whose inner vertices are
non-adjacent to $y_k$. Then we obtain a $1$-FAT obstruction. It remains to deal with the harder
situation where no such path exists.

Let $C(x_k)$ be the connected component of $G[Q] \setminus N[y_k]$ containing $x_k$. By our
assumption, $z_k \notin C(x_k)$. We denote by $W_k$ the subset of $N[y_k]$ containing those vertices
that are adjacent to some vertex of $C(x_k)$; see Fig.~\ref{fig:q_nodes_k_fat_statement}, middle.
Notice that the vertices of $C(x_k)$ appear only in sections and subtrees to the left of
$s_\beta(Q)$.  Therefore, every vertex of $W_k$ lies in the sections of $Q$ and stretches from the
left of $s_\beta(Q)$ to $s_\gamma(Q)$; see Fig.~\ref{fig:q_nodes_k_fat_statement}, right. In other
words, $W_k \subseteq s_\beta(Q) \cap s_\gamma(Q)$ and every vertex of $W_k$ is adjacent to $z_k$.

Let $C$ be a connected component of $G[Q] \setminus W_k$. If $C$ contains a vertex from some section
of $Q$, we call it \emph{big}. Notice that in this situation $C$ has a vertex contained in two
consecutive sections of $Q$ and their subtrees. Otherwise, $C$ consists of some vertices of a
subtree of $Q$, and we call it \emph{small}. The section above a subtree containing a small
component is a subset of $W_k$. Additionally, if two small components are placed in two different
subtrees, the two sections above these subtrees are different.

The graph $G[Q] \setminus W_k$ is disconnected, as $x_k$ and $z_k$ belong to different
components. Let us denote the connected component containing $y_k$ by $C(y_k)$, and the one
containing $z_k$ by $C(z_k)$.  
Let $t_k$ be a vertex of $W_k$ whose section
$s_{t_k}^\rt(Q)$ is leftmost. Let $P_k$ be a shortest path from $x_k$ to $t_k$ in $G[C(x_k) \cup \{t_k\}]$; see
Fig.~\ref{fig:q_nodes_k_fat_statement}, right. We distinguish two cases.

\begin{figure}[t!]
\centering
\includegraphics{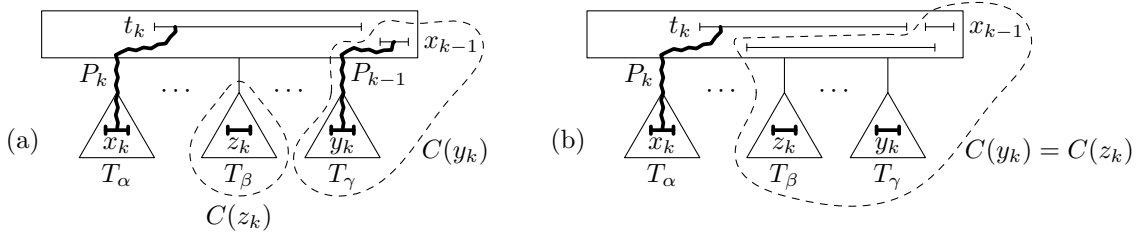}
\caption{(a) In Case 1, there exists a path $P_{k-1}$ from $x_{k-1}$ to $y_k$ whose inner
vertices avoid $z_k$. (b) In Case 2, we have
$C(y_k) = C(z_k)$ and such a path might no longer exist. For instance, every path from $x_{k-1}$ to
$y_k$ in $C(y_k)$ might use the depicted interval in the sections of $Q$, which is also
adjacent to $z_k$.}
\label{fig:q_nodes_k_fat_two_cases}
\end{figure}

\emph{Case 1: $C(y_k) \ne C(z_k)$.} This case is very similar to the proof of
Lemma~\ref{lem:q_node_two}; see Fig.~\ref{fig:q_nodes_k_fat_two_cases}a. Every vertex of $W_k$ is
adjacent to some vertex of $C(x_k)$ and to some vertex of $C(y_k)$. Therefore, it is also adjacent
to every vertex of $C(z_k)$. If $C(z_k)$ was big, then we could reverse its sections in the Q-node,
contradicting the fact that there are only two possible orderings for a Q-node. Therefore, $C(z_k)$
is small. Notice that then $C(y_k)$ is not small, since otherwise we would get $s_\beta(Q) = W_k =
s_\gamma(Q)$, contradicting Lemma~\ref{lem:different_sections}. Thus, $C(y_k)$ is big.

Let us set $y_{k-1} = z_k$ and $z_{k-1} = y_k$. The vertex $t_k$ is not universal for $C(y_k)$;
otherwise, every vertex of $W_k$ would be universal and this would give additional orderings of
$C(y_k)$ in $Q$.  Let $x_{k-1}$ be a vertex of $C(y_k) \setminus N[t_k]$ whose section
$s_{x_{k-1}}^\lft(Q)$ is leftmost. Notice that $s_{x_{k-1}}^\lft(Q)$ is always the next section to
$s_{t_k}^\rt(Q)$. Let $P_{k-1}$ be a shortest path from $x_{k-1}$ to $z_{k-1}$ in $C(y_k)$. By
Lemma~\ref{lem:monotone_paths}, all inner vertices of $P_{k-1}$ are adjacent to $t_k$.  Since this
path lies in $C(y_k)$, the inner vertices are non-adjacent to $y_{k-1}$, $x_k$ and $P_k$. We have
constructed a $2$-FAT obstruction.

\emph{Case 2: $C(y_k) = C(z_k)$.} In this case, the component $C(y_k)$ is big; see
Fig.~\ref{fig:q_nodes_k_fat_two_cases}b. Therefore, similarly as above, $t_k$ is not universal for
$C(y_k)$. We put $y_{k-1} = z_k$ and $z_{k-1} = y_k$. We choose $x_{k-1} \in C(y_k) \setminus
N[t_k]$ in the same way as in Case 1.  Notice that $x_{k-1}$ is a non-neighbor of $y_{k-1}$, since
otherwise it would be a neighbor of $t_k$.  On the other hand, $x_{k-1}$ might be adjacent to
$z_{k-1}$ or not. If it is, we get a $2$-FAT obstruction for $k=2$ with $P_1 = x_{k-1}z_{k-1}$. If
it is not, we proceed as follows.
	
As before, every shortest path from $x_{k-1}$ to $z_{k-1}$ has all inner vertices adjacent to $t_k$.
Since all vertices of $C(y_k)$ are non-adjacent to $x_k$ and the inner vertices of $P_k$, every
shortest path satisfies this as well.  There exists a shortest path from $x_{k-1}$ to
$z_{k-1}$ in $C(y_k)$, but we cannot guarantee that the inner vertices of this path are non-adjacent
to $y_{k-1}$. We solve this issue by applying the entire argument of the proof recursively to
$C(y_k)$.

In every representation extending the partial representation, the intervals of $C(x_k)$ form a
connected subset of the real line placed to the left of $\inter{y_k}'$. Therefore, $\inter{t_k}$
stretches from $C(x_k)$ to $\inter{z_k}'$, covering $\inter{y_k}'$. Thus $\inter{x_{k-1}}$ is placed
to the right of $\inter{z_k}' = \inter{y_{k-1}}'$ in every extending representation (see
Fig.~\ref{fig:k_FAT}b).  Again, $\inter{y_{k-1}}'$ has to be placed between $\inter{x_{k-1}}$ and
$\inter{z_{k-1}}'$. We assume that $\inter{x_{k-1}}$ is pre-drawn on the right of $\inter{y_{k-1}}'$
and repeat the same argument for $C(y_k)$ and the MPQ-tree restricted to these vertices. The role of
$x_k$, $y_k$ and $z_k$ is played by $x_{k-1}$, $y_{k-1}$ and $z_{k-1}$, respectively. (The ordering
of the pre-drawn intervals is flipped.)

The paragraphs above show the induction step of our proof (by induction on, say, the number of
considered sections of $Q$). By the induction hypothesis, we find a $(k-1)$-FAT obstruction.  By
making $x_{k-1}$ free and adding $x_k$, $t_k$ and $P_k$, we get a $k$-FAT obstruction in the
original partial representation.  Clearly $t_k$ is adjacent to the entire $(k-1)$-FAT obstruction
with the exception of $x_{k-1}$, since all further vertices are contained in a section to the left
of $s_{x_{k-1}}^\lft(Q)$. The reason is that we always use shortest paths which are $Q$-monotone by
Lemma~\ref{lem:monotone_paths}.  By the same reason, they are non-adjacent to the inner vertices of
$P_k$ and to $x_k$, as required.

To make the argument complete, we should check that all the assumptions used throughout the proof
apply recursively, in particular the arguments concerning non-universality of $t_{k-1}$ and
reversing big components. This is true because both components $C(y_{k-1})$ and $C(z_{k-1})$ of
$C(y_k) \setminus W_{k-1}$ appear to the left of $x_{k-1}$, so $t_k$ and the other vertices of $W_k$
are universal for them. This property is preserved throughout the recursion, so $C(y_\ell)$ and
$C(z_\ell)$ are adjacent to all vertices of $W_k,W_{k-1},\dots,W_{\ell+1}$. Similarly,
the rest of the inductive proof can be formalized.\qed
\end{proof}

\begin{figure}[t!]
\centering
\includegraphics{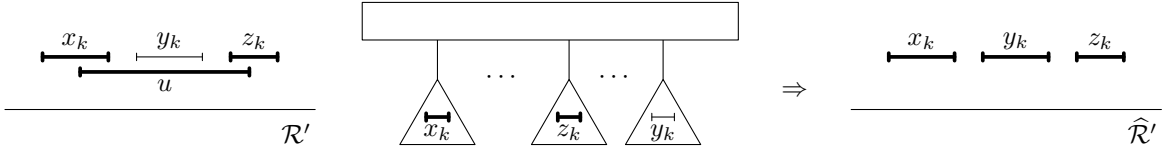}
\caption{Suppose that we show that a partial representation $\calR'$ has three pre-drawn intervals
as on the left, and that there is a vertex $y_k$ adjacent to $u$ and non-adjacent to $x_k$ and
$z_k$. Then $\inter{y_k}$ has to be placed between $\inter{x_k}'$ and $\inter{z_k}'$ in every
extending representation. Thus, we can assume it is pre-drawn there and obtain a modified partial
representation $\widehat\calR'$. If we further show that $x_k$, $y_k$ and $z_k$ are placed in
appropriate sections of $G[Q]$ for some Q-node $Q$, we can apply $k$-FAT Lemma~\ref{lem:k_fat} and
we get a $k$-FAT obstruction in $G[Q]$ and $\widehat\calR'[\{x_k,y_k,z_k\}]$. Together with $\inter
u'$, this forms a $k$-BI obstruction in $G$ and $\calR'$.}
\label{fig:q_nodes_using_k_fat}
\end{figure}

The above proof shows that the structure of a Q-node can be highly complicated, leading to 
complicated obstructions such as $k$-FAT. Actually, $k$-FAT Lemma~\ref{lem:k_fat} is a very useful tool
because it can be also applied in situations where not all
$x_k$, $y_k$, and $z_k$ are pre-drawn, to build other obstructions.
Fig.~\ref{fig:q_nodes_using_k_fat} shows an example.

\begin{lemma} \label{lem:k_fat_contains_1_fat}
Consider a $k$-FAT obstruction $H_k$ for $k>2$. If we swap the positions of $\inter{x_k}'$ and
$\inter{y_k}'$, then we obtain a new obstruction which contains a $1$-FAT obstruction for $x'_1 =
y_k$, $y'_1 = x_k$, and $z'_1 = z_k$. Further, if $k=2$ and this does not happen, then $x_2$ is
adjacent to $t_2$.
\end{lemma}

\begin{proof}
For $k \ge 3$, the graph $H_k \setminus N[x_k]$ is connected; in particular, there
exists a path $P'_1 = y_kt_{k-1}z_k$ avoiding $N[x_k]$. For $k=2$, there exists a path $P'_1 =
y_2t_2z_2$ avoiding $N[x_2]$, unless $x_2$ is adjacent to $t_2$.\qed
\end{proof}

\begin{figure}[b!]
\centering
\includegraphics{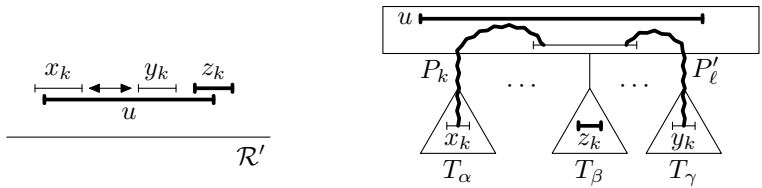}
\caption{When $\inter u'$ single overlaps $\inter{z_k}'$, and the vertices $x_k$, $y_k$, and $z_k$ are
placed in the MPQ-tree as on the right, we get a $(k,\ell)$-CE obstruction.}
\label{fig:q_nodes_kl_ce_statement}
\end{figure}

\heading{$\boldsymbol{(k,\ell)}$-CE Lemma.} Suppose that we have the situation in
Fig.~\ref{fig:q_nodes_kl_ce_statement}. We can easily show that there is some $(k,\ell)$-CE
obstruction by applying $k$-FAT Lemma~\ref{lem:k_fat} twice, once when $\inter{x_k}$ is on the left
of $\inter{y_k}$ and once when it is on the right.  The following lemma reveals its structure in
detail.

\begin{lemma}[$(\boldsymbol{k},\boldsymbol{\ell})$-CE] \label{lem:kl_ce}
Let $Q$ be a Q-node with children $T_1,\dots,T_n$, and let $a$, $b$ and $c$ be three cliques of
$T[Q]$ contained respectively in $T_\alpha$, $T_\beta$ and $T_\gamma$, for $\alpha < \beta < \gamma$.
Let $x_k \in a$, $y_k \in c$ and $z_k \in P(b)$ be three non-adjacent vertices having a common
pre-drawn neighbor $u$ such that $\inter{u}'$ single overlaps $\inter{z_k}'$.
Then $G[Q] \cup \{u\}$ and $\calR'[\{z_k,u\}]$ contain a $(k,\ell)$-CE obstruction, where either $\ell=1$ or
$k = \ell = 2$.
\end{lemma}

\begin{proof}
The simplest case is when there exist a path $P_k$ from $x_k$ to $z_k$ avoiding $N[y_k]$, and a path
$P_\ell'$ from $y_k$ to $z_k$ avoiding $N[x_k]$. Let $P_k$ and $P_\ell'$ be shortest such paths as
in Fig.~\ref{fig:q_nodes_kl_ce_statement}, right. We get a $(1,1)$-CE obstruction. By
Lemma~\ref{lem:monotone_paths}, the paths $P_k$ and $P_\ell'$ are monotone. Therefore, their inner
vertices are non-adjacent to each other, with the possible exception of the last vertices before
$z_k$, which can be adjacent or even identical. Concerning minimality, we can always find one of the
three finite $(1,1)$-CE obstructions depicted in Fig.~\ref{fig:kl_CE_obstructions}a. The reason is
that when paths $P_k$ and $P'_\ell$ are long, we can take as $x_k$ and $y_k$ one of their inner
vertices, making them shorter.

Suppose next that there exists no path $P_k$ from $x_k$ to $z_k$ avoiding $N[y_k]$.  Let $C(x_k)$,
$W_k$, and $t_k$ be defined as in the proof of $k$-FAT Lemma~\ref{lem:k_fat}. Following the argument
in that proof, we get the subgraph $H_k$ of a $k$-FAT obstruction, which is not the complete $k$-FAT
obstruction because $x_k$ and $y_k$ are free. 

\emph{Case 1: There exists some path $P_\ell'$ from $y_k$ to $z_k$ avoiding $N[x_k]$.} Let $P'_\ell$
be a shortest such path (notice that $\ell = 1$). Together with the above subgraph $H_k$, we get a $(k,1)$-CE
obstruction; see Fig.~\ref{fig:q_nodes_kl_ce_two_cases}a. In particular, if some vertex $w \in W_k$
is non-adjacent to $x_k$ (possibly $w = t_k$), we can use $P'_\ell = y_kwz_k$. 

We note that when $k \ge 3$, such a path $P_\ell'$ necessarily exists, as we can use $P'_\ell =
y_kt_{k-1}z_k$, as argued in Lemma~\ref{lem:k_fat_contains_1_fat}. Therefore, the $(k,1)$-CE
obstructions consist of the subgraph $H_k$ together with $u$; assuming minimality, we have that
either $u$ is adjacent to all vertices of $H_k$, or $u = t_k$.  If $k = 2$, then $P_\ell'$ might
still exist but it might be longer and might use inner vertices not contained in $H_k$. Concerning
minimality, we always find one of the three $(2,1)$-CE obstructions depicted in
Fig.~\ref{fig:kl_CE_obstructions}b. Indeed, $P_2$ can be assumed to be of length one or two, since
otherwise we could use one of its inner vertices as $x_2$. For length two, we get $P'_1 =
y_kt_kz_k$. For length one, we get a path $P'_1$ from $z_2$ to $y_2$, and we can assume that $y_2$
is adjacent to $x_1$ (otherwise we could use as $y_2$ the neighbor of $x_1$ on $P_1$).

\begin{figure}[t!]
\centering
\includegraphics{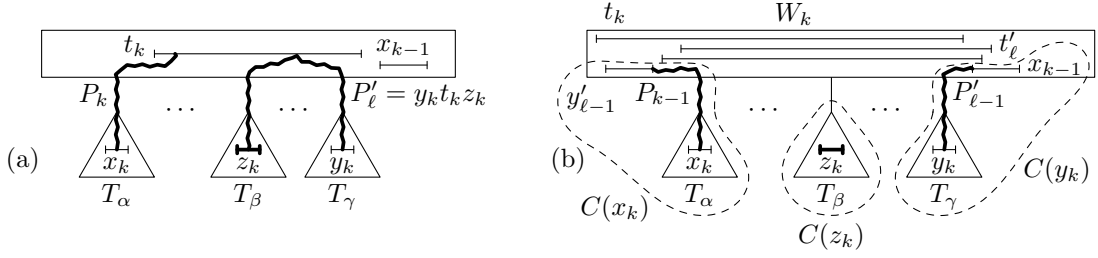}
\caption{(a) Case 1: If there exists a path $P'_\ell$ from $y_k$ to $z_k$ avoiding $N[x_k]$, then we get a
$(k,1)$-CE obstruction. (b) Case 2: We get a $(2,2)$-CE obstruction.}
\label{fig:q_nodes_kl_ce_two_cases}
\end{figure}

\emph{Case 2: No such path $P'_\ell$ exists.}
By Lemma~\ref{lem:k_fat_contains_1_fat}, $k = 2$.
We want to show that there exists a $(2,2)$-CE obstruction. 

Notice that all vertices $w\in W_k$ are
adjacent to $x_k$, $y_k$, and $z_k$, since otherwise there would exist a path $P'_\ell = y_kwz_k$.
Hence the vertices of $W_k$ belong to sections of $Q$, covering all subtrees between $T_\alpha$ and
$T_\gamma$; see Fig.~\ref{fig:q_nodes_kl_ce_two_cases}b. Let $C(y_k)$ and $C(z_k)$ be the components
of $G[Q] \setminus W_k$ containing $y_k$ and $z_k$, respectively.  Since there exists no path
$P'_\ell$, we obtain that $C(x_k)$, $C(y_k)$, and $C(z_k)$ are pairwise different.  To determine the
structure of a $(2,2)$-CE obstruction, we apply the argument from Case 1 of the proof of
$k$-FAT Lemma~\ref{lem:k_fat} symmetrically twice.

Let $t_k$ be a vertex of $W_k$ having leftmost section $s_{t_k}^\rt(Q)$ and let $t'_\ell$ be a
vertex of $W_k$ having rightmost section $s_{t'_\ell}^\lft(Q)$ (possibly $t_k = t'_\ell$).  It is
easy to see that $C(z_k)$ is small, otherwise we could flip it and obtain an
ordering of the maximal cliques not compatible with the Q-node.

Similarly as in the proof of $k$-FAT Lemma~\ref{lem:k_fat}, this implies that both $C(x_k)$ and
$C(y_k)$ are big. Therefore, $t_k$ is not universal for $C(y_k)$ and $t'_\ell$ is not universal for
$C(x_k)$. As in the proof of $k$-FAT Lemma~\ref{lem:k_fat}, we choose $x_{k-1} \in C(y_k)$
non-adjacent to $t_k$  and $y'_{\ell-1} \in C(x_k)$ non-adjacent to $t'_\ell$.  There exist paths
$P_{k-1}$ from $x_{k-1}$ to $y_k$ and $P'_{\ell-1}$ from $y'_{\ell-1}$ to $x_k$.  In consequence, we
obtain a $(2,2)$-CE obstruction.

Regarding minimality, notice that we can assume that $y_2$ is adjacent to $x_1$, and $x_2$ is adjacent to $y'_1$;
otherwise, we could choose as $y_2$ and $x_2$ the neighbors of $x_1$ and $y'_1$ on the paths $P_1$ and
$P'_1$, respectively. We get the four minimal finite $(2,2)$-CE obstructions that are illustrated in
Fig.~\ref{fig:kl_CE_obstructions}c.
\qed
\end{proof}

\begin{figure}[t!]
\centering
\includegraphics{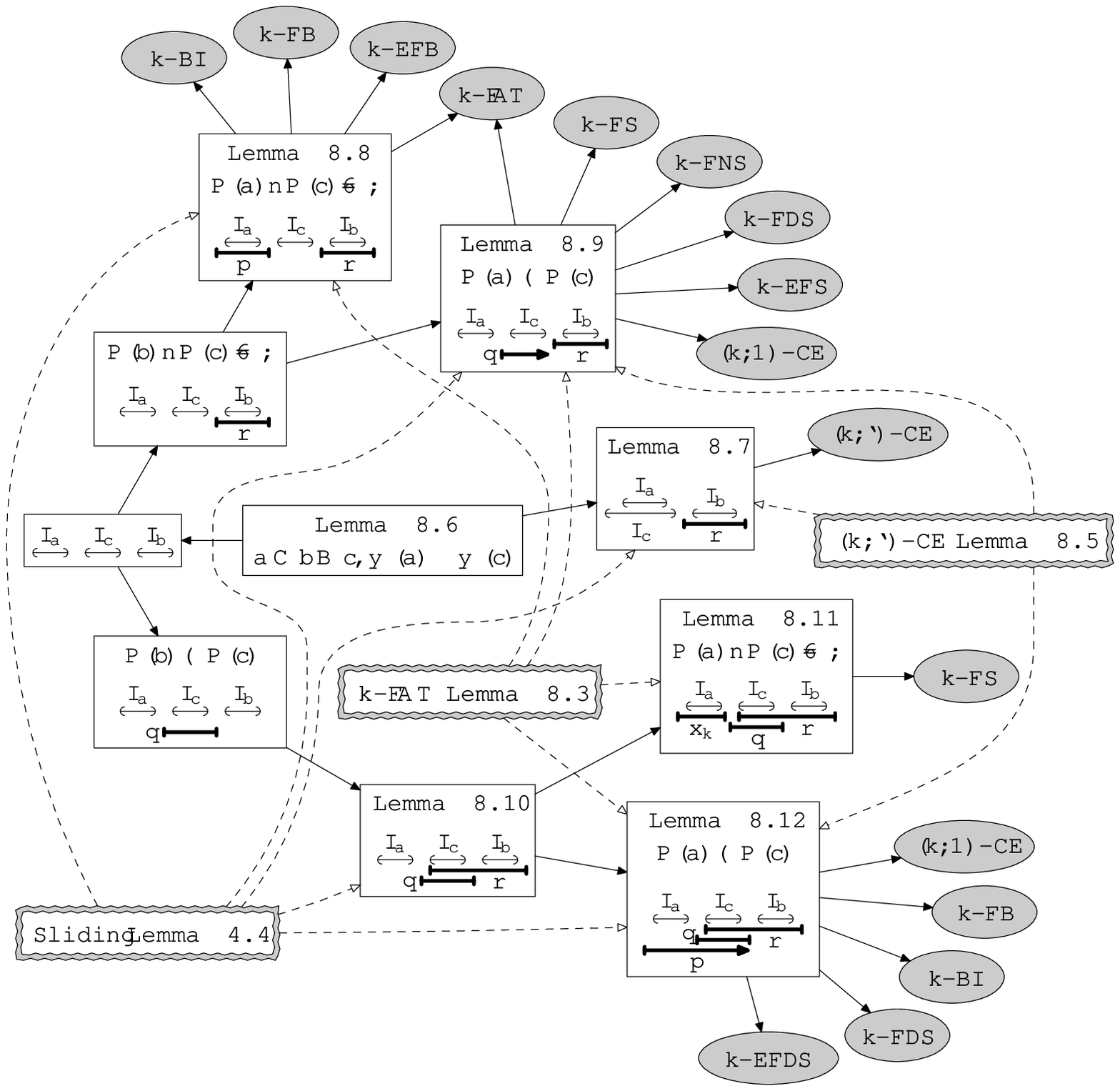}
\caption{A summary of Section~\ref{sec:q_nodes_three_subtrees}. The diagram starts in the middle with
Lemma~\ref{lem:q_node_wlog_lemma}. Inside the cases, we draw the positions of $I_a$, $I_b$, $I_c$,
and some pre-drawn intervals. An arrow at a pre-drawn interval means that it may be further
stretched in the given direction.  The obtained obstructions are highlighted in gray, the used tools
have highlighted borders.}
\label{fig:q_nodes_three_subtrees_diagram}
\end{figure}

\subsection{Cliques in Three Different Subtrees} \label{sec:q_nodes_three_subtrees}

When a Q-node $Q$ is obstructed by three maximal cliques $a \in T_\alpha$, $b \in T_\beta$ and
$T_\gamma$, where $\alpha < \beta < \gamma$, the situation is quite complex.
Fig.~\ref{fig:q_nodes_three_subtrees_diagram} gives an overview of the cases and obstructions
obtained in this case. 

\begin{lemma} \label{lem:q_node_wlog_lemma}
Without loss of generality, we can assume that $a \wlt b \wgt c$ and $\clqr(a) \le \clqr(c)$.
\end{lemma}

\begin{proof}
Since $a$, $b$ and $c$ create an obstruction, $b$ is either a minimal or a maximal element in $\wlt
|_{\{a,b,c\}}$. Without loss of generality (using the flip operation), we can assume that $b$ is
maximal, so $a \wlt b \wgt c$.  Since we can swap $a$ and $c$ by reversing the Q-node, we can assume
that $\clqr(a) \le
\clqr(c)$.\qed
\end{proof}

Since $a \wlt b \wgt c$, both $I_a$ and $I_c$ appear on the left of $I_b$.  Since $\clqr(a) \le
\clqr(c)$, either $I_a$ contained in $I_c$, or $I_a$ is on the left of $I_c$. The first case is easier:

\begin{lemma} \label{lem:q_node_containment_case}
If $I_a$ is contained in $I_c$, then $G$ and $\calR'$ contain a $(k,\ell)$-CE obstruction, where
$\ell = 1$ or $2 \ge k \ge \ell$.
\end{lemma}

\begin{proof}
The proof is illustrated in Fig.~\ref{fig:q_nodes_containment_case}.  By
Lemma~\ref{lem:open_intervals}, $P(c) \subseteq P(a)$. Since $b$ is placed between $a$ and $c$ in
the Q-node $Q$, every vertex contained in both $a$ and $c$ is contained in $b$ as well. Hence $P(c)
\subsetneq P(b)$, and there exists $r \in P(b) \setminus P(c)$. Since $\inter r'$ is on the right of
$I_c$, it is also on the right of $I_a$, and thus $r \notin P(a)$.

Let $u \in P^{\endlft}(c)$. We apply Sliding Lemma~\ref{lem:sliding} to $I_c$, $I_b$, and $\inter
r'$. We get a pre-drawn interval $\inter{z_k}'$ covering $r(u)$, and an induced path $P_{r,z_k}$
from $r$ to $z_k$ consisting of pre-drawn intervals not in $P(c)$. Therefore $z_k$ is on the left of
$c$ in $Q$. Since all pre-drawn intervals of $P_{r,z_k}$ do not belong to $P(c)$, they are on the
right of $I_c$. Thus they are also on the right of $I_a$, which implies that they do not belong to
$P(a)$. Consequently, $z_k$ is between $a$ and $c$ in $Q$.

\begin{figure}[t!]
\centering
\includegraphics{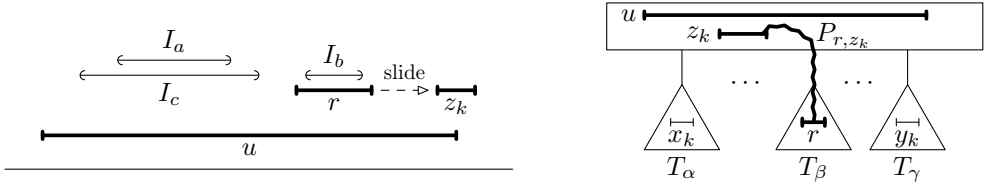}
\caption{Proof of Lemma~\ref{lem:q_node_containment_case}. On the left, the pre-drawn
intervals. On the right, their positions in the MPQ-tree.}
\label{fig:q_nodes_containment_case}
\end{figure}

Let $x_k \in a$ and $y_k \in c$ be vertices non-adjacent to $z_k$. By $(k,\ell)$-CE
Lemma~\ref{lem:kl_ce}, $x_k$, $y_k$, $z_k$, and $u$ create a $(k,\ell)$-CE obstruction, for $\ell =
1$ or $2 \ge k \ge \ell$. Notice that the clique associated to $z_k$ is some $b'\neq b$.\qed
\end{proof}

The case where $I_a$ is on the left of $I_c$ is further divided into several subcases. In the next two lemmas, we focus on the situation where $P(b) \setminus P(c) \ne \emptyset$.

\begin{lemma} \label{lem:q_node_p_and_r_case}
If $I_a$ is on the left of $I_c$, $P(b) \setminus P(c) \ne \emptyset$ and $P(a) \setminus P(c) \ne
\emptyset$, then $G$ and $\calR'$ contain a $k$-FAT, $k$-BI ($k \le 2$), $k$-FB, or $k$-EFB
obstruction.
\end{lemma}

\begin{proof}
The proof is illustrated in Fig.~\ref{fig:q_nodes_p_and_r_case}.  Let $p \in P(a) \setminus P(c)$ and
$r \in P(b) \setminus P(c)$. Then $\inter p'$ is on the left of $I_c$, and $\inter r'$ is on the
right of $I_c$.  Clearly, $\inter p'$ and $\inter r'$ are disjoint, so $p$ appears in the Q-node on
 the left of $r$. Let $u \in P^{\endlft}(c)$ and $v \in P^{\startrt}(c)$ (possibly $u = v$).

\begin{figure}[b!]
\centering
\includegraphics{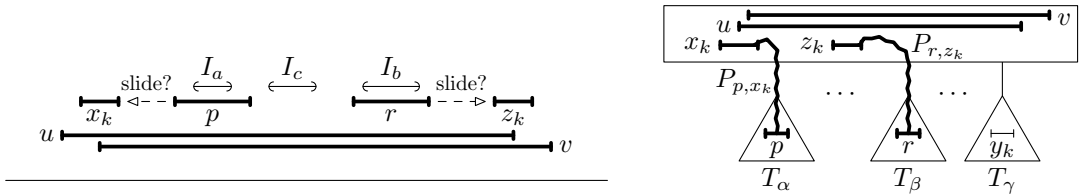}
\caption{Proof of Lemma~\ref{lem:q_node_p_and_r_case}. On the left, the pre-drawn
intervals, with possible sliding on each side. On the right, their positions in the MPQ-tree.}
\label{fig:q_nodes_p_and_r_case}
\end{figure}

If $r(u) \le r(r)$, then $z_k = r$. Obviously, $z_k$ is between $p$ and $c$ in the Q-node.
Otherwise, $P(c) \subsetneq P(b)$, and we apply Sliding Lemma~\ref{lem:sliding} to $I_c$, $I_b$, and
$r$. We obtain a pre-drawn interval $\inter{z_k}'$ not contained in $P(c)$ covering $r(u)$, and a
path $P_{r,z_k}$ whose inner vertices are pre-drawn and not contained in $P(c)$. Notice that all of
these pre-drawn vertices are on the right of $I_c$.  Therefore, these vertices are not contained in
$P(a)$. In this case, we also get that $z_k$ is between $p$ and $c$ in the Q-node.  

Similarly, if $\ell(p) \le \ell(v)$, then $x_k = p$. Otherwise, we use the flipped version of
Sliding Lemma~\ref{lem:sliding} to $I_c$, $I_a$, and $p$, which gives a pre-drawn interval
$\inter{x_k}'$ not contained in $P(c)$ covering $\ell(v)$. By a similar argument, in both cases, we
show that $x_k$ is on the left of $z_k$ in the Q-node.

Let $y_k \in c$ be a vertex non-adjacent to $z_k$ (possibly, $y_k = u$ or $y_k = v$).  Such a vertex
exists because $z_k$ is on the left of $c$ in the Q-node. Notice that $y_k$ is also non-adjacent to
$x_k$. Since $y_k$ is adjacent to $u$ and $v$, in every extending representation $\inter{y_k}$ is
between $\inter{x_k}'$ and $\inter{z_k}'$. So we can assume that it is pre-drawn in this position
and, by $k$-FAT Lemma~\ref{lem:k_fat}, we get a $k$-FAT obstruction. Together with $u$ and $v$ (or
possibly only one of them), we get one of the obstructions in
Fig.~\ref{fig:q_nodes_p_and_r_case_obstructions}.\qed

\begin{figure}[t!]
\centering
\includegraphics{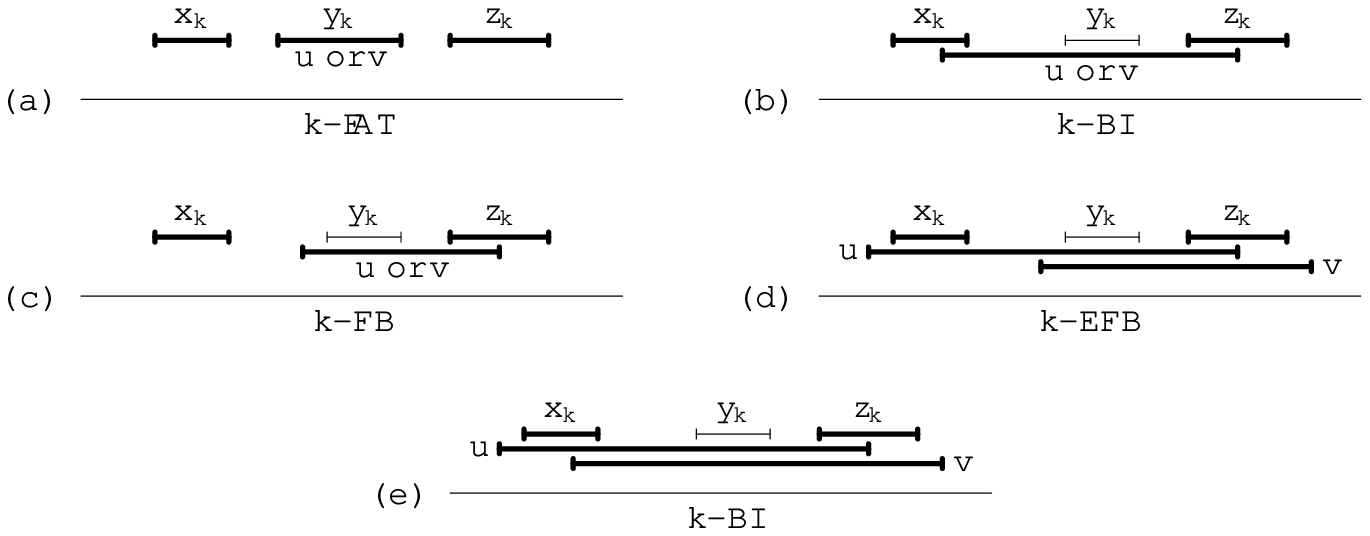}
\caption{The different obstructions obtained in the proof of Lemma~\ref{lem:q_node_p_and_r_case}.
If $\ell(x_k) \le \ell(u) \le r(u) \le r(z_k)$, we get one of the obstructions (a) to (c). Since
$z_k$ is in the Q-node between $x_k$ and $c$, if $u$ or $v$ intersect $x_k$, then they also
intersect $z_k$. In the cases (d) and (e), $\ell(u) < \ell(x_k)$ and $r(z_k) < r(v)$. Since $u$
intersects $z_k$, there are only two possible obstructions.}
\label{fig:q_nodes_p_and_r_case_obstructions}
\end{figure}

\end{proof}

\begin{lemma} \label{lem:q_node_q_and_r_case}
If $I_a$ is on the left of $I_c$, $P(b) \setminus P(c) \ne \emptyset$ and $P(a) \subsetneq
P(c)$, then $G$ and $\calR'$ contain a $k$-FAT, $(k,1)$-CE, $k$-FS, $k$-FDS, $k$-FNS, or $k$-EFS obstruction.
\end{lemma}

\begin{proof}
We choose $r \in P(b) \setminus P(c)$ and $q \in P(c) \setminus P(a)$ with leftmost right endpoint.
Then $\inter r'$ is on the right of $I_c$ and $\inter q'$ is on the right of $I_a$.  We note that
$q$ might be adjacent to $r$ or not, and might belong to $P(b)$ or not. Since $P(a) \subsetneq
P(c)$, we get from the structure of the Q-node that also $P(a) \subsetneq P(b)$.  Let $u \in
P^{\endlft}(a)$. Notice that at least one of $q$ and $u$ belongs to $P^{\endlft}(c)$. 

\emph{Case 1: $u \in P^{\endlft}(c)$}. Then $r(u) \le r(q)$ and $P(c) \subsetneq P(b)$; the
situation is depicted in Fig.~\ref{fig:q_nodes_q_and_r_case}a.  We apply Sliding
Lemma~\ref{lem:sliding} to $I_c$, $I_b$, and $r$. We get a pre-drawn interval $z_k \notin P(c)$
covering $r(u)$, and a path $P_{r,z_k}$ consisting of pre-drawn intervals not contained in $P(c)$.
Therefore, $z_k$ is on left of $c$ in the Q-node. Since $I_a$ is on the left of $I_c$, all vertices
of $P_{r,z_k}$ are also not contained in $P(a)$. Thus $z_k$ is on the right of $a$ in the Q-node.

Choose $y_k \in c$ non-adjacent to $z_k$.  By Lemma~\ref{lem:non_adjacent}, there exists $x_k \in a$
non-adjacent to both $z_k$ and $q$.  Since $z_k$ is between $a$ and $c$ in the Q-node, also $x_k$ is
non-adjacent to $y_k$.  Since $y_kqz_k$ is a path avoiding $N[x_k]$, by $(k,\ell)$-CE
Lemma~\ref{lem:kl_ce} we obtain a $(k,1)$-CE obstruction.

\begin{figure}[b!]
\centering
\includegraphics{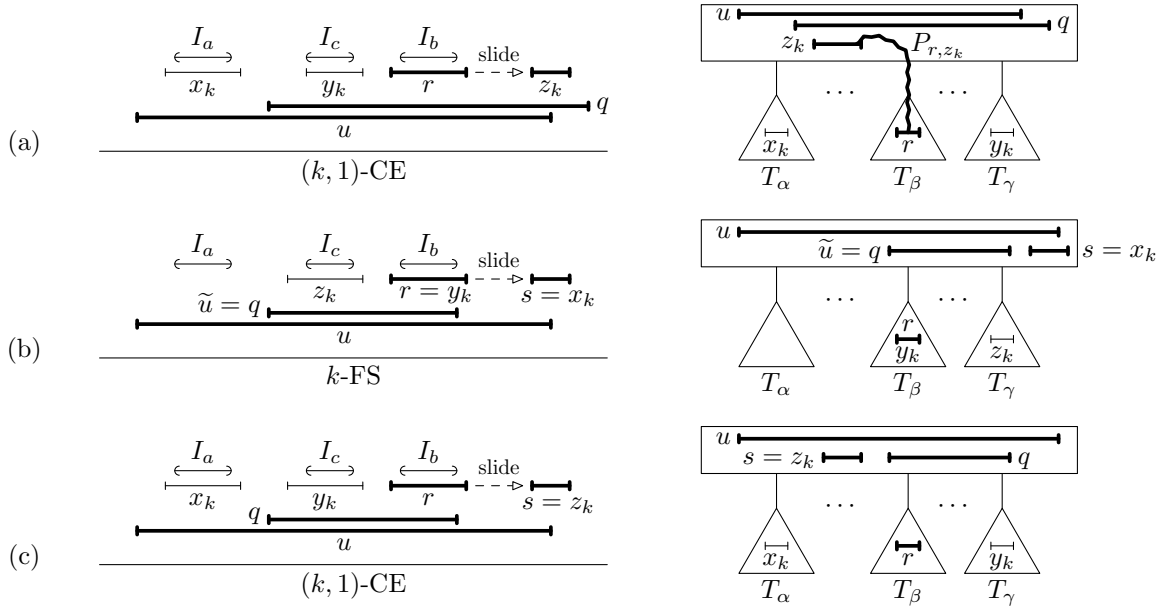}
\caption{Proof of Lemma~\ref{lem:q_node_q_and_r_case}. On the left, the pre-drawn
intervals. On the right, their positions in the MPQ-tree. (a) Case 1. (b) Subcase 2A. (c) Subcase 2B.}
\label{fig:q_nodes_q_and_r_case}
\end{figure}

\emph{Case 2: $q \in P^{\endlft}(c)$}. Then $r(q) < r(u)$.
First we argue that, without loss of generality, we can assume that either $\inter q'$ and $\inter
r'$ are disjoint, or $\inter r'$ covers $r(q)$. Suppose that $\inter r'$ is contained in $\inter
q'$. Since $q \in P^{\endlft}(c)$, this implies that $P(c) \subsetneq P(b)$. By applying Sliding
Lemma~\ref{lem:sliding} to $I_c$, $I_b$ and $r$, we obtain a pre-drawn interval $\widetilde r$ not
in $P(c)$ which covers $r(q)$.  We also get a path $P_{r,\widetilde r}$ whose vertices are
pre-drawn and not contained in $P(c)$; since $I_a$ is on the left of $I_c$, they are also not in $P(a)$.
Therefore $\widetilde r$ is between $a$ and $c$ in the Q-node. Further, $\widetilde r$ belongs to
some clique $\widetilde b$ for which $I_{\widetilde b}$ is on the right of $I_c$. From now on, we
work with $\widetilde r$ as $r$, and with $\widetilde b$ as $b$. Hence the
assumption on the relative positions of $\inter q'$ and $\inter r'$ holds.

We apply Sliding Lemma~\ref{lem:sliding} to $I_a$, $I_b$ and $r$, and we get a pre-drawn interval $s
\notin P(a)$ covering $r(u)$. This sliding is weaker that in Case 1: we know that $s$ is on the
right of $a$, but we do not know its position with respect to $c$.  We distinguish three subcases
according to the relative positions of $q$ and $s$ in the Q-node.

\emph{Subcase 2A: $s$ is on the right of $q$.}
The situation is depicted in Fig.~\ref{fig:q_nodes_q_and_r_case}b.  Let $x_k = s$ and $y_k = r$.  If
$\inter q'$ is on the left of $\inter r'$, let $z_k = q$. Otherwise, let $z_k \in c$ be a vertex
non-adjacent to $r$, but possibly adjacent to $x_k$.  Since in every extending
representation $z_k$ is placed on the left of $y_k$, we can apply $k$-FAT
Lemma~\ref{lem:k_fat} to $x_k$, $y_k$ and $z_k$, and get a subgraph $H_k$.  If $\inter q'$ is on the
left of $\inter r'$, then $H_k$ gives a $k$-FAT obstruction.  If $\inter r'$ covers $r(q)$, then
$H_k$ together with $\widetilde u = q$ gives a $k$-FS obstruction.

\begin{figure}[t!]
\centering
\includegraphics{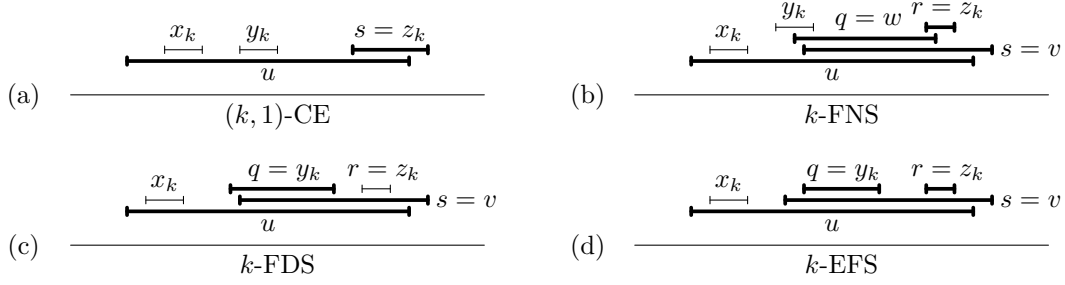}
\caption{Four possible obstructions obtained in Subcase 2C of the proof of
Lemma~\ref{lem:q_node_q_and_r_case}. (a) If $s \notin P(c)$, we get a $(k,1)$-CE obstruction. (b) If
$\inter q'$ intersects $\inter r'$, we get a $k$-FNS obstruction. Recall that the relative order of
$\ell(q)$ and $\ell(s)$ does not matter. (c) If $\inter q'$ is on the left of $\inter r'$ and
$\ell(q) \le \ell(s)$, we get a $k$-FDS obstruction.  (d) If $\inter q'$ is on the left of $\inter
r'$ and $\ell(s) < \ell(q)$, we get a $k$-EFS obstruction.}
\label{fig:q_nodes_q_and_r_case_obstructions}
\end{figure}

\emph{Subcase 2B: $s$ is on the left of $q$.}
We choose $x_k \in a$ and $y_k \in c$ non-adjacent to $s$; such vertices exist because $s$ is
between $a$ and $c$ in the Q-node. By $(k,\ell)$-CE Lemma~\ref{lem:kl_ce}, we get a $(k,\ell)$-CE
obstruction for $x_k$, $y_k$, $z_k = s$ and $u$. Notice that we can construct a path
$P_{y_k,z_k}$ from $y_k$ to $z_k$ avoiding $N[x_k]$ by applying Sliding Lemma~\ref{lem:sliding} to
$I_a$, $I_c$, and $q$. Thus $\ell = 1$.

\emph{Subcase 2C: $\inter s'$ intersects $\inter q'$.}
Notice that $\inter s'$ also intersects $\inter r'$. Therefore, if $s \notin P(c)$, then it appears
in the Q-node between $a$ and $c$. Let $z_k = s$, we get a $(k,1)$-CE obstruction as follows.  We
choose $y_k \in c$ non-adjacent to $z_k$.  By Lemma~\ref{lem:non_adjacent}, there exists $x_k \in a$
non-adjacent to $q$, $y_k$, and $z_k$. By $(k,\ell)$-CE Lemma~\ref{lem:kl_ce}, we get a $(k,1)$-CE
obstruction for $x_k$, $y_k$, $z_k$ and $u$ as illustrated in
Fig.~\ref{fig:q_nodes_q_and_r_case_obstructions}a; notice that the path $y_kqz_k$ avoids $N[x_k]$.

It remains to deal with the situation when $s \in P(c)$.  Let $z_k = r$. If $\inter q'$ intersects
$\inter {r}'$, let $y_k \in c$ be a vertex non-adjacent to $r$; otherwise let $y_k = q$. By
Lemma~\ref{lem:non_adjacent}, there exists $x_k \in a$ non-adjacent to $q$, $y_k$, and $z_k$.  In
every extending representation, $\inter{y_k}$ is placed on the left of $\inter{z_k}'$, and
$\inter{x_k}$ is placed on the left of $\inter{y_k}$.  Therefore, by $k$-FAT Lemma~\ref{lem:k_fat},
we get a subgraph $H_k$ of a $k$-FAT obstruction.  Together with $u$, $v = s$, $w = q$ (for $y_k \ne
q$), or possibly some of them, we get a $k$-FDS, $k$-EFS, or $k$-FNS obstruction; see
Fig.~\ref{fig:q_nodes_q_and_r_case_obstructions}b, c, and d.\qed
\end{proof}

The case where $P(b) \subsetneq P(c)$ is addressed in Lemmas~\ref{lem:q_node_p_and_q_case}
and~\ref{lem:q_node_q_case}. First, we need an auxiliary result.

\begin{lemma} \label{lem:q_node_q_case_wlog}
If $I_a$ is on the left of $I_c$ and $P(b) \subsetneq P(c)$, there exist $q \in P(c) \setminus
P(b)$ and $r \in P(b) \setminus P(a)$ such that $\inter q'$ is on the right of $I_a$ and on the left
of $I_b$, and $\inter r'$ is on the right of $I_a$, containing $I_c$ and $I_b$.  Without loss of
generality, $\inter q'$ covers $\ell(r)$.
\end{lemma}

\begin{proof}
The proof is depicted in Fig.~\ref{fig:q_nodes_q_case_wlog}. Clearly, there exists $q \in P(c) \setminus
P(b)$. Due to the structure of the Q-node, we also have that $q \notin P(a)$.  Therefore, $\inter
q'$ is between $I_a$ and $I_b$. 

\begin{figure}[b!]
\centering
\includegraphics{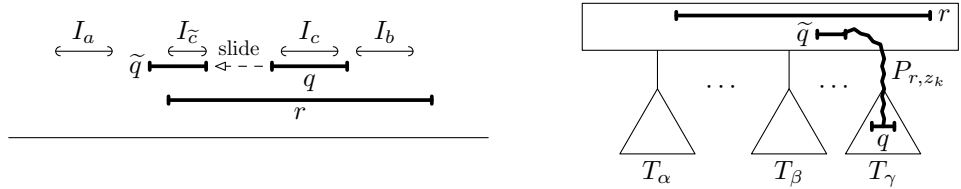}
\caption{Proof of Lemma~\ref{lem:q_node_q_case_wlog}. On the left, the pre-drawn
intervals. On the right, their positions in the MPQ-tree.}
\label{fig:q_nodes_q_case_wlog}
\end{figure}

Next, we argue that there exists $r \in P(b) \setminus P(a)$. For contradiction, assume that $P(b)
\subsetneq P(a)$. Let $v \in P^{\startrt}(b)$; notice that $v$ contains $I_a$ and $I_c$. By the
flipped version of Sliding Lemma~\ref{lem:sliding} applied to $I_c$, $I_b$ and $q$, there exists a
path consisting of pre-drawn intervals not contained in $P(b)$ from $q$ to $z$, where $\inter z'$
covers $\ell(v)$. At least one interval of this path intersects $I_a$, so it belongs to $P(a)$. This
contradicts the fact that $b$ is between $a$ and $c$ in the Q-node. Hence, there exists $r \in P(b)
\setminus P(a)$.

We choose $r$ having rightmost left endpoint. Clearly, $\inter r'$ is on the right of $I_a$, and
contains $I_b$ and $I_c$. Suppose that $\ell(r) < \ell(q)$. Since $r$ has  rightmost left endpoint
among all intervals in $P(b) \setminus P(a)$, and every interval in $P(a)$ has its left endpoint
more to the left, we obtain that $r \in P^{\startrt}(b)$. Therefore, we can apply the flipped
version of Sliding Lemma~\ref{lem:sliding} to $I_c$, $I_b$ and $q$. We get a pre-drawn interval
$\widetilde q \notin P(b)$ covering $\ell(r)$, and a path $P_{q,\widetilde q}$ from $q$ to
$\widetilde q$ whose vertices are not in $b$.  Therefore, $\widetilde q$ is on the right of $b$ in
the Q-node. Let $\widetilde c$ be a maximal clique containing $\widetilde q$. Since $I_{\widetilde
c}$ is contained in $\widetilde q$, it is between $I_a$ and $I_b$. Therefore, we can work with
$\widetilde q$ and $\widetilde c$ instead of $q$ and $c$. Thus we can assume that $\inter q'$ covers
$\ell(r)$.\qed
\end{proof}

For $P(b) \subsetneq P(c)$, we distinguish two cases.

\begin{lemma} \label{lem:q_node_p_and_q_case}
If $I_a$ is on the left of $I_c$, $P(b) \subsetneq P(c)$, and $P(a) \setminus P(c) \ne
\emptyset$, then $G$ and $\calR'$ contain a $k$-FS obstruction.
\end{lemma}

\begin{proof}
The proof is illustrated in Fig.~\ref{fig:q_nodes_p_and_q_case}. By Lemma~\ref{lem:q_node_q_case_wlog}, there exist $q \in P(c) \setminus
P(b)$ and $r \in P(b) \setminus P(a)$ such that $\inter q'$ covers $\ell(r)$. Let $x_k \in P(a) \setminus P(c)$
and $y_k = q$. Then $\inter{x_k}'$ is on the left of $I_c$, and therefore also on the left of $I_b$. Thus
$x_k \notin P(b)$. We infer that $x_k$ is on the left of $b$ in the Q-node, so it is non-adjacent to
$y_k$. In consequence, $\inter{x_k}'$ is on the left of $\inter{y_k}'$. Let $z_k \in b$ be a vertex
non-adjacent to $y_k$. 

\begin{figure}[t!]
\centering
\includegraphics{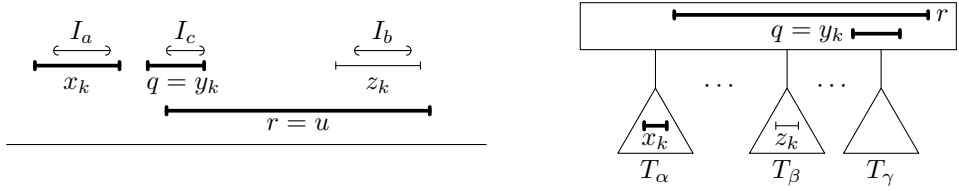}
\caption{Proof of Lemma~\ref{lem:q_node_p_and_q_case}. We derive that $\inter{x_k}'$ is on the
left of $\inter q'$, which gives a $k$-FS obstruction based on the positions in the Q-node.}
\label{fig:q_nodes_p_and_q_case}
\end{figure}

If $z_k$ is adjacent to $x_k$, we get a $1$-FS obstruction.  Otherwise, in every extending
representation, $\inter{z_k}$ is to the right of $\inter{y_k}'$.  By $k$-FAT
Lemma~\ref{lem:k_fat}, we get a subgraph $H_k$ of a $k$-FAT obstruction.  Together with $u = r$,
this leads to a $k$-FS obstruction.\qed
\end{proof}

\begin{lemma} \label{lem:q_node_q_case}
If $I_a$ is on the left of $I_c$, $P(b) \subsetneq P(c)$ and $P(a) \subsetneq P(c)$, then $G$ and $\calR'$ contain a
$(k,1)$-CE, $k$-FB, $k$-BI, $k$-FDS, or $k$-EFDS obstruction.
\end{lemma}

\begin{proof}
Let $p \in P^{\endlft}(a)$, and $q$ be the vertex from Lemma~\ref{lem:q_node_q_case_wlog}. By
applying Sliding Lemma~\ref{lem:sliding} to $I_a$, $I_c$ and $q$, we get a pre-drawn interval $s
\notin P(a)$ covering $r(p)$, and path $P_{q,s}$ of intervals not in $P(a)$, so $s$ appears on the
right of $a$ in the Q-node.  Similarly, as in Case 2 of the proof of
Lemma~\ref{lem:q_node_q_and_r_case}, we distinguish three cases according to the relative positions
of $s$ and $q$ in the Q-node; see Fig~\ref{fig:q_nodes_q_case}.

\begin{figure}[t!]
\centering
\includegraphics{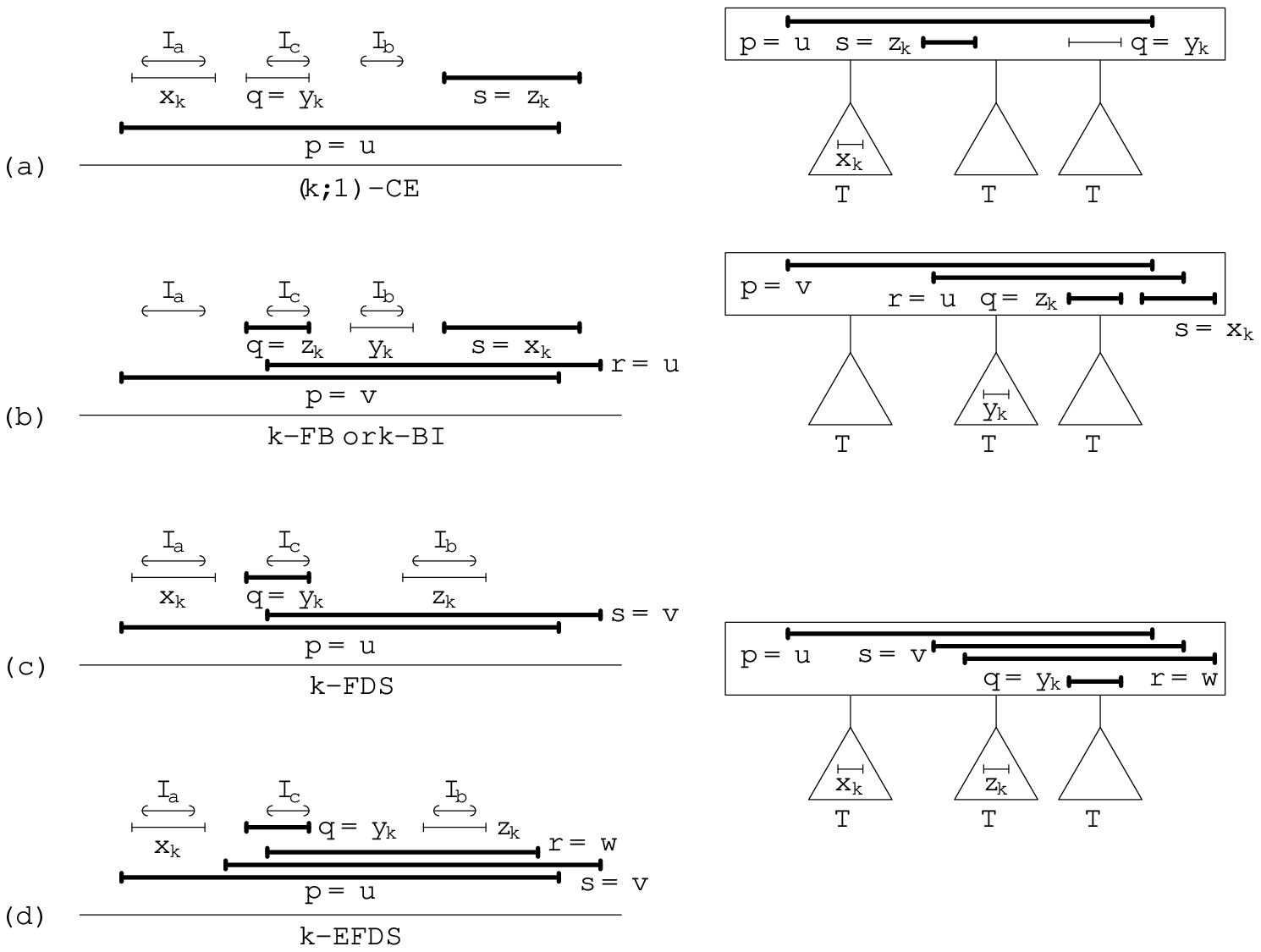}
\caption{Proof of Lemma~\ref{lem:q_node_q_case}. On the left, the pre-drawn
intervals. On the right, their positions in the MPQ-tree. (a) Case 1. (b) Case 2. (c) Case 3, if
$\ell(y_k) \le \ell(s)$. (d) Case 3, if $\ell(y_k) > \ell(s)$.}
\label{fig:q_nodes_q_case}
\end{figure}

\emph{Case 1: $s$ is on the left of $q$.} By Lemma~\ref{lem:non_adjacent}, there exists $x_k \in a$
non-adjacent to all vertices of $P_{q,s}$, in particular non-adjacent to $s$ and $q$.  Let $y_k =
q$, $z_k = s$, and $u = p$.  Clearly, $z_k$ is between $x_k$ and $y_k$ in the Q-node. By
$(k,\ell)$-CE Lemma~\ref{lem:kl_ce} and the existence of $P_{y_k,z_k}$, we get a $(k,1)$-CE
obstruction. Notice that $\inter{y_k}'$ can be made free; see Fig.~\ref{fig:q_nodes_q_case}a.

\emph{Case 2: $s$ is on the right of $q$.} Since $q$ is between $b$ and $s$ in the Q-node, we get
that $s \notin P(b)$. Let $x_k = s$, $z_k = q$, and $u = r$, where $r$ is the vertex from
Lemma~\ref{lem:q_node_q_case_wlog}. There exists $y_k \in b$ non-adjacent to $z_k$ and, by the
structure of the Q-node, also non-adjacent to $x_k$.  Since $y_k$ is adjacent to $p$ and $r$,
$\inter{y_k}$ is between $\inter{x_k}'$ and $\inter{z_k}'$ in every extending representation.  By
$k$-FAT Lemma~\ref{lem:k_fat}, we get a subgraph $H_k$ of a $k$-FAT obstruction.  If $r(r) \le
r(x_k)$, together with $u$, we obtain a $k$-FB or a $k$-BI obstruction. If $r(r) > r(x_k)$, together
with $u$ and $v=p$, we obtain a $k$-BI obstruction; see Fig.~\ref{fig:q_nodes_q_case}b.

\emph{Case 3: $\inter s'$ intersects $\inter q'$.}
Since $s$ contains $I_b$, it belongs to $P(b)$. Let $y_k = q$, $u = p$, $v = s$, and $w = r$; we
note that possibly $s=r$.  By Lemma~\ref{lem:non_adjacent}, there exists $x_k \in a$ non-adjacent to
$y_k$, $v$, and $w$. Since $x_k$ is adjacent to $u$, then $\inter{x_k}$ is on the left of
$\inter{y_k}'$ in every extending representation. Finally, there exists $z_k \in b$ non-adjacent to
$y_k$. Since $z_k$ is adjacent to $u$, $v$, and $w$, we have that $\inter{z_k}$ is on the right of
$\inter{y_k}'$ in every extending representation.

Since $z_k$ is between $x_k$ and $y_k$ in the Q-node, we can apply $k$-FAT Lemma~\ref{lem:k_fat},
which gives a subgraph $H_k$ of a $k$-FAT obstruction.  If $\ell(y_k) \le \ell(v)$, together with
$u$ and $v$, we obtain a $k$-FDS obstruction; see Fig.~\ref{fig:q_nodes_q_case}c. If $\ell(y_k) >
\ell(v)$, together with $u$, $v$, and $w$, we get a $k$-EFDS obstruction; see
Fig.~\ref{fig:q_nodes_q_case}d.\qed
\end{proof}

In summary, we conclude:

\begin{lemma}[The Q-node, Three Subtrees] \label{lem:q_node_three}
If the three cliques creating the obstruction belong to three different subtrees, then $G$ and
$\calR'$ contain a $k$-FAT, $k$-BI ($k \le 2$), $k$-FS, $k$-EFS, $k$-FB, $k$-EFB, $k$-FDS,
$k$-EFDS, $k$-FNS, or $(k,\ell)$-CE obstruction (either $k=\ell=2$, or $k \ge \ell = 1$).
\end{lemma}

\begin{proof}
For an overview, see the diagram in Fig.~\ref{fig:q_nodes_three_subtrees_diagram}. The proof follows
from Lemmas~\ref{lem:q_node_wlog_lemma}, \ref{lem:q_node_containment_case},
\ref{lem:q_node_p_and_r_case}, \ref{lem:q_node_q_and_r_case}, \ref{lem:q_node_q_case_wlog},
\ref{lem:q_node_p_and_q_case}, and \ref{lem:q_node_q_case}.\qed
\end{proof}

%% file: main_result.tex
\section{Proofs of the Main Results} \label{sec:main_result}

Now, we are ready to put all results together to prove the main theorem. It states that a partial
representation $\calR'$ of $G$ is extendible if and only if $G$ and $\calR'$ contain none of the
obstructions described in Section~\ref{sec:definition_of_obstructions}.

\begin{proof}[Theorem~\ref{thm:repext_obstructions}]
If $G$ and $\calR'$ contain one of the obstructions, they are non-extendible by
Lemma~\ref{lem:correctness_all}. It remains to prove the converse. If $G$ is not an interval graph,
it contains an LB obstruction~\cite{lb_graphs}. Otherwise, $G$ is an interval graph and there exists
an MPQ-tree $T$ for it. By Lemma~\ref{lem:repext_char}, we know that a partial representation
$\calR'$ is extendible if and only if $T$ can be reordered according to $\wlt$. If it cannot be
reordered, then the reordering algorithm fails in some node of $T$. If this reordering fails in a
leaf, we get a $1$-BI obstruction by Lemma~\ref{lem:leaf}. If it fails in a P-node, we get an SE,
$1$-BI, or $1$-FAT obstruction by Lemma~\ref{lem:p_node}. And if it fails in a Q-node, we get one of
the obstructions of Section~\ref{sec:definition_of_obstructions} by
Lemmas~\ref{lem:q_node_three_cliques}, \ref{lem:q_node_two}, and \ref{lem:q_node_three}.\qed
\end{proof}

Next, we show that a partial representation $\calR'$ is extendible if and only if every quadruple of
pre-drawn intervals is extendible by itself.

\begin{proof}[Corollary~\ref{cor:helly_obstructions}]
The result follows from the fact that all the obstructions of Theorem~\ref{thm:repext_obstructions}
contain at most four pre-drawn intervals.\qed
\end{proof}

Concerning the certifying algorithm, we first show that $k$-FAT obstructions can be found in linear
time:

\begin{lemma} \label{lem:computing_k_fat}
Suppose that the assumptions of $k$-FAT Lemma~\ref{lem:k_fat} are satisfied. Then we can find a
$k$-FAT obstruction in time $\O(n+m)$.
\end{lemma}

\begin{proof}
Since the proof of $k$-FAT Lemma~\ref{lem:k_fat} is constructive, the algorithm follows it.  Let $Q$
be the Q-node.  We search the graph $G[Q] \setminus N[y_k]$ from $x_k$ to compute $C(x_k)$, and test
whether $z_k$ belongs to it. If it does, the algorithm stops and outputs $1$-FAT. Otherwise, we compute
$W_k$, choose $t_k$, and store it together with $P_k$.
We choose $x_{k-1}$ as in the the proof; if $s_i(Q) = s_{t_k}^\rt(Q)$, then either
$s_{x_{k-1}}^\lft(Q) = s_{i+1}(Q)$, or $x_{k-1}$ belongs to sections of $T_{i+1}$.
Then we apply the rest of the algorithm recursively.  It is important that then we can remove
$C(x_k)$ and $W_k$ from the graph because they are not used in the remainder of the obstruction.

Since the algorithm searches each vertex and edge of $G[Q] \setminus N[y_k]$ at most once when
computing $C(x_j)$, we obtain that the algorithm runs in time $\O(n+m)$.\qed
\end{proof}

Similarly, a $(k,\ell)$-CE obstruction can be obtained from $(k,\ell)$-CE Lemma~\ref{lem:kl_ce} in
time $\O(n+m)$.  Since obstructions are built constructively, we get a linear-time certifying
algorithm for the partial representation extension problem:

\begin{proof}[Corollary~\ref{cor:certifying_algorithm}]
We can assume that $G$ is an interval graph; otherwise we can find an LB obstruction in time
$\O(n+m)$ using~\cite{finding_lb_graphs}.  Each interval graph has $\O(n)$ maximal cliques of total
size $\O(n+m)$, and that they can be found in linear time~\cite{recog_chordal_graphs}.  We compute
the MPQ-tree $T$ in time $\O(n+m)$ using~\cite{korte_mohring}. Next, we use the partial
representation extension algorithm of~\cite{kkosv} in time $\O(n+m)$, which either finds an
extending representation, or finds an obstructed node which cannot be reordered according to $\wlt$.
We distinguish three cases according to the distinct types of obstructed nodes.

\emph{Case 1: A leaf cannot be reordered.}
We output a $1$-BI obstruction in time $\O(n)$, by searching the partial representation.

\emph{Case 2: A P-node $P$ cannot be reordered.}
From the partial representation extension algorithm, we get directly a two-cycle, ensured by
Lemma~\ref{lem:p_node_two_cycle}, and four maximal cliques $a$, $b$, $c$, and $d$ defining it. By
Lemma~\ref{lem:p_node_three_cliques}, one of these maximal cliques can be omitted, and
it can be clearly found in constant time. It remains to output an SE, $1$-BI, or $1$-FAT
obstruction in time $\O(n+m)$, by following Lemma~\ref{lem:p_node}. For $1$-BI and $1$-FAT
obstructions, we find a shortest path in $G[P] \setminus N[y_k]$ by searching the graph.

\emph{Case 3: A Q-node $Q$ cannot be reordered.}
From the partial representation extension algorithm, we get four maximal cliques defining the
obstruction and, by following Lemma~\ref{lem:q_node_three_cliques}, we can reduce it to at most three
maximal cliques. An SE obstruction can be computed in time $\O(n+m)$. If three maximal
cliques are contained in two subtrees, we follow Lemma~\ref{lem:q_node_two} and output one of the
obstructions in time $\O(n+m)$.

If three maximal cliques belong to three different subtrees, we follow the structure of the proof of
Lemma~\ref{lem:q_node_three}. In all cases, we derive some vertices somehow placed in the Q-node and
some pre-drawn intervals, which can be easily done in time $\O(n+m)$. Next, we either apply $k$-FAT
Lemma, or $(k,\ell)$-CE Lemma to construct the obstruction, which can be done in time $\O(n+m)$ by
Lemma~\ref{lem:computing_k_fat}.\qed
\end{proof}

%% file: conclusions.tex
\section{Conclusions} \label{sec:conclusions}

In this paper, we have described the minimal obstructions that make a partial interval representation
non-extendible. There are three main points following from the proof:

\begin{packed_enum}
\item Minimal obstructions for the partial representation extension problem are much more
complicated than minimal forbidden induced subgraphs of interval graphs, characterized by
Lekkerkerker and Boland~\cite{lb_graphs}.
\item Nevertheless, it is possible to describe these obstructions using structural results derived
in~\cite{kkosv} and in this paper. We show that almost all of these obstructions consist of three
intervals $x_k$, $y_k$ and $z_k$ that are forced by the partial representation to be drawn in an
incorrect left-to-right order. This incorrect placement leads to the complex zig-zag structure of a
$k$-FAT obstruction. 
\item The structure of the sections of a Q-node $Q$ can be very intricate. Suppose that we contract
in $G[Q]$ the sections of each subtree $T_i$ into one vertex. Then we get an interval graph which
has a unique interval representation up to flipping the real line. Such interval graphs have been
extensively studied, see for instance~\cite{hanlon,fishburn_unique,prescribed_lengths_npc}.
Therefore, our structural results needed to find minimal obstructions may be of independent
interest.
\end{packed_enum}

\heading{Structural Open Problems.}
The first open problem we propose is a characterization of minimal obstructions for other graph
classes. We select those classes for which polynomial-time algorithms are
known~\cite{cfk,kkkw,kkorssv}.  \emph{Circle graphs} (\circle) are intersection graphs of chords of
a circle. \emph{Function graphs} (\fun) are intersection graphs of continuous functions $f : [0,1]
\to \mathbb R$, and \emph{permutation graphs} (\perm) are function graphs which can be represented
by linear functions.  \emph{Proper interval graphs} (\pint) are intersection graphs of closed
intervals in which no interval is a proper subset of another interval. \emph{Unit interval graphs}
(\uint) are intersection graphs of closed intervals of length one. 

\begin{problem} \label{prob:other_minimal_obstructions}
What are the minimal obstructions for partial representation extension of the classes \circle, \fun,
\perm, \pint, and \uint?
\end{problem}

The second open problem involves a generalization of partial representations called \emph{bounded
representations}~\cite{bko,kkorssv,soulignac}. Suppose that a graph $G$ is given together with two closed
intervals $L_v$ and $R_v$ for every vertex $v \in V(G)$. A \emph{bounded representation} of $G$ is a
representation such that $\ell(v) \in L_v$ and $r(v) \in R_v$ for every vertex $v \in V(G)$. We call
bounds \emph{solvable} if and only if there exists a bounded representation.
This generalizes partial representations: we can use singletons $L_v$ and $R_v$ for pre-drawn 
intervals and put them equal $\mathbb R$ for the others.

\begin{problem} \label{prob:bounded_representations}
What are minimal obstructions making bounds for interval graphs unsolvable?
\end{problem}

\heading{Algorithmic Open Problems.}
We have described a linear-time certifying algorithm that can find one of the minimal obstructions
in a non-extendible partial representation. There are several related computational problems,
suggested by Jan Kratochv\'{\i}l, for which the complexity is open:

\begin{problem} \label{prob:obstruction_testing}
What is the computational complexity of the problem of testing whether a given minimal obstruction
is contained in $G$ and $\calR'$?
\end{problem}

Since a minimal obstruction contains at most four pre-drawn intervals, we can test over all subsets
of at most four pre-drawn intervals whether they form an obstruction (say, by freeing the rest of
them and testing whether the modified partial representation is extendible). If $k$ is fixed, we can
test whether the subgraph of a given obstruction is contained in $G$. Given a triple $x_k$, $y_k$
and $z_k$ forming a $k$-FAT obstruction, the proof of k-FAT Lemma~\ref{lem:k_fat} and the algorithm
of Lemma~\ref{lem:computing_k_fat} constructs it while minimizing $k$. The approach needs to be
changed to check whether they also form an $\ell$-FAT obstruction, for $\ell > k$.

The next problem generalizes the partial representation extension problem.

\begin{problem} \label{prob:representation_modification}
What is the computational complexity of testing whether at most $\ell$ pre-drawn intervals can be
freed to make a partial representation extendible $\calR'$?
\end{problem}

Similar problems are usually \cNP-complete. On the other hand, we propose the following reformulation
which might lead to a polynomial-time algorithm. Every minimal obstruction contains at most four
pre-drawn intervals. Let $P$ be the set of pre-drawn intervals, and let $\calS$ consist of all
subsets of $P$ of size at most four which form an obstruction.  We can clearly compute $\calS$ in
polynomial time. Then the problem above is equivalent to finding a minimal hitting set of $P$ and
$\calS$. This problem is in general \cNP-complete, but the extra structure given by the MPQ-tree
might make it tractable.

\begin{problem} \label{prob:graph_modification}
What is the complexity of testing whether it is possible to remove at most $\ell$ vertices from an
interval graph $G$ to make a partial representation extendible $\calR'$?
\end{problem}

This problem is fundamentally different from Problem~\ref{prob:representation_modification}, in
which the partial representation $\calR'$ is modified. In this problem, we modify the graph $G$
itself, changing its structure. When we remove a pre-drawn vertex, we also remove its pre-drawn
interval from the partial representation. We note that the assumption that $G$ is an interval graph
is important. For general graphs $G$, the problem is known to be \cNP-complete even when $\calR' =
\emptyset$~\cite{node_deletion_hereditary}.